\documentclass{pspum-l}
\usepackage{amscd}  
\usepackage[all]{xy}

\newtheorem{theorem}{Theorem}[section]
\newtheorem{proposition}[theorem]{Proposition}
\newtheorem{corollary}[theorem]{Corollary}
\newtheorem{lemma}[theorem]{Lemma}
\newtheorem{conjecture}[theorem]{Conjecture}

\theoremstyle{definition}
\newtheorem{definition}[theorem]{Definition}
\newtheorem{example}[theorem]{Example}

\theoremstyle{remark}
\newtheorem{remark}[theorem]{Remark}

\numberwithin{equation}{section}

\begin{document}

\title{The Chow Ring of a Classifying Space}
\author{Burt Totaro}
\address{Department of Mathematics, University of Chicago,
5734 S.~University Ave., Chicago, IL 60637}
\email{totaro@math.uchicago.edu}
\thanks{The author was supported in part by a Sloan Research Fellowship.}

\date{November 24, 1997}



\def\Z{\mbox{\bf Z}}
\def\Q{\mbox{\bf Q}}
\def\R{\mbox{\bf R}}
\def\C{\mbox{\bf C}}
\def\Proj{\mbox{\bf P}}
\def\arrow{\rightarrow}
\def\inj{\hookrightarrow}
\def\imp{\Rightarrow}
\def\im{\mbox{im}}
\def\k{\overline{k}}
\def\QQ{\overline{\Q}}
\def\qed{\ QED. \vspace{3mm}}
\def\invlim{\varprojlim}
\def\surj{\rightarrow}
\def\BPr{BP\langle r \rangle}
\def\gr{\text{gr}}
\def\Rep{\text{Rep}}
\def\sphat{^{\wedge}}
\def\ch{\text{ch}}

\maketitle

For any linear algebraic group $G$,
we define a ring $CH^*BG$, the ring of
characteristic classes with values in the Chow ring 
(that is, the ring of
algebraic cycles modulo rational equivalence) for principal
$G$-bundles over smooth algebraic varieties. We show that
this coincides with the Chow ring of any quotient variety $(V-S)/G$
in a suitable range of dimensions, where $V$ is a representation of $G$
and $S$ is a closed subset such that $G$ acts freely outside $S$.
As a result, computing the Chow ring of $BG$ amounts to the computation
of Chow groups for a natural class of algebraic varieties with a lot
of torsion in their cohomology. Almost nothing is known about this
in general. For $G$ an algebraic group over the complex numbers, there is
an obvious ring homomorphism $CH^*BG\arrow H^*(BG,\Z )$. Less obviously,
using the results of \cite{TotaroMU},
this homomorphism factors through the quotient of the complex cobordism ring
$MU^*BG$ by the ideal generated by the elements of negative degree
 in the coefficient ring
$MU^*=\Z[x_1,x_2,\dotsc]$,
that is, through the ring $MU^*BG\otimes_{MU^*}\Z$. 
(For clarity, let
us mention that the classifying space of a complex algebraic group
is homotopy equivalent to that of its maximal compact subgroup. Moreover,
every compact Lie group arises as the maximal compact subgroup of a
unique complex reductive group.)

The most interesting
result of this paper is that in all the examples where we can compute
the Chow ring of $BG$, it maps isomorphically to the topologically defined
ring $MU^*BG\otimes_{MU^*}\Z$. Namely, this is true for
finite abelian groups, the symmetric groups, tori, $GL(n,\C)$, $Sp(2n,\C)$,
$O(n)$, $SO(2n+1)$, and $SO(4)$. (The computation of the Chow ring
for $SO(2n+1)$ is the result of discussion between me and
Rahul Pandharipande.
Pandharipande then proceeded to compute the Chow ring of $SO(4)$,
which is noticeably more difficult \cite{Pandharipande}.) We also
get various additional information about the Chow rings of the symmetric
groups, the group $G_2$, and so on.

Unfortunately, the map $CH^*BG\arrow MU^*BG\otimes_{MU^*}\Z$
is probably not always an isomorphism, since this would imply in particular
that $MU^*BG$ is concentrated in even degrees. By Ravenel, Wilson,
and Yagita,
$MU^*BG$ is concentrated in even degrees if all the Morava $K$-theories
of $BG$ are concentrated in even degrees \cite{RWY}. The latter statement,
for all compact Lie groups $G$,
was a plausible conjecture of Hopkins, Kuhn, and
Ravenel \cite{HKR}, generalizing the theorem of
Atiyah-Hirzebruch-Segal on the topological $K$-theory of $BG$ \cite{Atiyah},
\cite{AS}.
But it has now been disproved by Kriz \cite{Kriz}, using
the group $G$ of strictly
upper triangular $4 \times 4$ matrices over $\Z/3$.

Nonetheless, there are some reasonable conjectures to make. First, if
$G$ is a complex algebraic group such that the complex cobordism ring
of $G$ is concentrated in even degrees, say after tensoring with
$\Z_{(p)}$ ($\Z$ localized at $p$) for a fixed prime number $p$,
then the homomorphism
$$CH^*BG\arrow MU^*BG\otimes_{MU^*}\Z$$
should become an isomorphism after tensoring with $\Z_{(p)}$.
Second, we can hope that the Chow ring of $BG$ for any complex algebraic
group $G$ has the good properties which complex cobordism was formerly
expected to have. In particular, the Chow ring of $BG$ for a finite group $G$
should be additively generated
by transfers to $G$ of Chern classes of representations of subgroups
of $G$ (see section \ref{surjcases}). This conjecture suggests thinking
of the Chow ring of $BG$ as a better-behaved substitute for the
ordinary group cohomology ring.

In many cases where we compute the Chow ring of $BG$, it maps injectively
to the integral cohomology of $BG$: that is in fact true for all
the examples mentioned above (finite abelian groups, the symmetric groups,
tori, $GL(n,\C)$, $Sp(2n,\C)$, $O(n)$, $SO(2n+1)$, and $SO(4)$).
 In such cases, we can say that the Chow ring
of $BG$ is the image of $MU^*BG$ in $H^*(BG,\Z)$. Nonetheless, that is not
true in general: for $G=\Z/2 \times SO(4)$, the Chow ring of $BG$ is equal
to $MU^*BG\otimes_{MU^*}\Z$, but this does not inject into $H^*(BG,\Z)$.
This observation was used to give a topological construction of
elements of the Griffiths group of certain smooth projective
varieties in \cite{TotaroMU}.

The motivation for looking at these questions is the problem of 
determining which torsion cohomology classes on a smooth projective
variety $X$ can be represented by algebraic cycles. In particular,
the image of $CH^1X\arrow H^2(X,\Z )$ always contains the torsion
subgroup of $H^2(X,\Z )$, but Atiyah and Hirzebruch
\cite{Atiyah-Hirzebruch} showed that for
$i\geq 2$ there are varieties $X$ and torsion elements of $H^{2i}(X,\Z )$
which are not in the image of $CH^iX$. Their examples are Godeaux-Serre
varieties, that is, quotients of smooth complete intersections by free
actions of finite groups $G$. This leads to the question of actually computing
the Chow groups of Godeaux-Serre varieties to the extent possible.
It turns out that this problem almost completely reduces to the problem
of computing the Chow ring of $BG$ if we believe a suitable version
of a conjecture of Nori's (section \ref{godeaux}).
For the codimension 2 Chow group of Godeaux-Serre varieties,
we can prove a strong relation to the codimension 2 Chow group
of $BG$.

The outline of the paper is as follows. In section \ref{def}, we
define the Chow ring of $BG$ for an algebraic group $G$ over
any field. In section \ref{general}, restricting to groups over the complex
numbers, we construct the factorization $CH^*BG\arrow MU^*BG\otimes_{MU^*}\Z
\arrow H^*(BG,\Z)$. The geometric work has already been done in
\cite{TotaroMU}, so unfortunately what remains is a technical
argument about inverse limits. Section \ref{K_0} describes the analogous
notion of the algebraic $K$-group $K_0$ of $BG$, which has been completely
computed by Merkurjev \cite{Merkurjev}.
 We use this to show, for example, that the image of
the homomorphism $CH^iBG\arrow 
(MU^*BG\otimes_{MU^*}\Z)^{2i}$ contains $(i-1)!$ times the latter group.
Moreover, the map is an isomorphism when $i$ is 1 or 2.
Section \ref{surjcases} discusses the conjecture that $CH^*BG$
is generated by transferred Euler classes for $G$ finite.
In section \ref{godeaux}, we conjecture the relation of the Chow ring
of $BG$ to the Chow ring of Godeaux-Serre varieties, and we prove most of it
in the case of $CH^2$. Section \ref{Kunneth} discusses the Chow groups
of products in the cases we need. Sections \ref{CPintro} to \ref{trans}
lead up to the proof that $CH^*BS_n\arrow MU^*BS_n\otimes_{MU^*}\Z$
is an isomorphism, where $S_n$ denotes the symmetric group. Section
\ref{mod2} gives a curious description of the ring $CH^*BS_n\otimes\Z/2$, and
section \ref{fields} analyzes the Chow ring of the symmetric group
over base fields other than the complex numbers. Section \ref{bound}
proves the last general result on the Chow ring of 
classifying spaces, an explicit
 upper bound
for the degree of generators for this ring. Nothing similar is known
for either ordinary cohomology or complex cobordism. The bound is used
in sections \ref{classical} and \ref{other} to compute $CH^*BG$ for most
of the classical groups $G$.

The inspiration for the definition of the Chow ring of $BG$
came from Bogomolov's work on the problem of rationality
for quotient varieties \cite{Bogomolov}.
I thank Dan Edidin, H\'{e}l\`{e}ne Esnault, Bill Fulton, Bill Graham,
Bruno Kahn, Nick Kuhn, Doug Ravenel, Frank Sottile,
Steve Wilson, and Nobuaki Yagita
for telling me about  their work.
I have included Chris Stretch's unpublished
proof that $MU^*BS _6$ is not generated by Chern classes
in section \ref{surjcases}.

\tableofcontents

\section{Definition of the Chow ring of an algebraic group}
\label{def}

Fix a field $k$. We will work throughout in the category of 
separated schemes of finite type over $k$, called ``varieties'' for short.
(We will try to say explicitly when we use the word ``variety'' in the
more standard sense of a reduced irreducible separated scheme of finite
type over $k$.)
For any variety $X$, let $CH_iX$ denote
the group of $i$-dimensional
algebraic cycles on $X$ modulo rational equivalence \cite{Fulton},
called the $i$th Chow group of $X$. If $X$ is smooth of dimension $n$,
we define the $i$th Chow cohomology group of $X$ to be
$CH^iX:=CH_{n-i}X$. Then $CH^*X$ is a graded ring, called the Chow ring
of $X$.
The Chow ring is contravariant for arbitrary maps of smooth varieties.
Let $G$ be a linear algebraic group over $k$. 
A principal $G$-bundle over an algebraic variety
$X$ is a variety $E$ with free $G$-action such  that $X=E/G$.

\begin{theorem}
\label{indep}
Let $G$ be a linear algebraic group over a field $k$.
Let $V$ be any representation of $G$ over $k$
such that $G$ acts freely
outside a $G$-invariant
closed subset $S\subset V$ of codimension $\geq s$. 
Suppose that the geometric quotient $(V-S)/G$ (in the sense of
\cite{Mumford}) exists as a
variety over $k$. Then the ring $CH^*(V-S)/G$, restricted to degrees
less than $s$, is independent (in a canonical way) of the representation $V$
and the closed subset $S$.
\end{theorem}

Moreover, such pairs $(V,S)$ exist with the codimension of $S$ in $V$
arbitrarily
 large (see Remark \ref{quot}, below). So we can make the following definition:

\begin{definition}
For a linear algebraic group $G$ over a field $k$, define 
$CH^iBG$ to be the group 
$CH^i(V-S)/G$ for any $(V,S)$ as in Theorem \ref{indep}
such that $S$ has codimension greater than $i$ in $V$.
\end{definition}

Then $CH^*BG$ forms a ring by Theorem \ref{indep}, called the Chow
ring of $BG$. The ring $CH^*BG$ depends on the field $k$; when we want
to indicate this explicitly in the notation, we will write $CH^*(BG)_k$.

 For reductive groups $G$, such as finite groups,
 we can provide a more convincing
justification for the name: the Chow ring of $BG$ is equal to the ring
of characteristic classes for principal $G$-bundles over smooth
quasi-projective varieties, in the following sense.

\begin{theorem}
\label{twodefs}
Let $G$ be a reductive group over a field $k$. Then the above group
$CH^iBG$ is naturally identified  with the set of assignments $\alpha$
to every smooth quasi-projective variety $X$ over $k$
with a principal $G$-bundle
$E$ over $X$ of an element $\alpha(E)\in CH^iX$, such that for any map
$f:Y\arrow X$ we have $\alpha(f^*E)=f^*(\alpha (E))$. The ring structure
on $CH^*BG$ is the obvious one.
\end{theorem}

\begin{remark}
\label{quot}
If $G$ is finite, then the geometric quotient $V/G$ exists
as an affine variety for all representations $V$ of $G$ \cite{Mumford}.
So $(V-S)/G$ is a quasi-projective variety for all closed subsets $S\subset
V$ such that $G$ acts freely on $V-S$.

For any  linear algebraic group $G$ and any positive integer $s$, there
is a representation $V$ of $G$ and a closed subset $S\subset V$
of codimension at least $s$ such that $G$ acts freely on $V-S$ and
$(V-S)/G$ exists as a quasi-projective variety over $k$. For example,
let $W$ be any faithful representation of $G$, say of dimension $n$,
 and let $V=\text{Hom}
(A^{N+n},W)\cong W^{\oplus N+n}$ for $N$ large. Let $S$ be the closed subset
in $V$ of non-surjective linear maps $A^{N+n}\arrow W$.
 Then the codimension of
$S$ in $V$ goes to infinity as $N$ goes to infinity. Also, $(V-S)/G$
exists as a quasi-projective variety. Indeed, we can view it as the homogeneous
space $GL(N+n)/(H\times G)$ where $H$ is the subgroup of $GL(N+n)$
which acts as the identity on a given $n$-dimensional quotient space of 
$A^{N+n}$, so that $H$ is an extension of $GL(N)$ by an additive group
$\text{Hom}(A^N,A^n)$. Here any quotient of a linear algebraic group
by a closed subgroup exists as a quasi-projective variety
by \cite{Waterhouse},  pp.~121-122. Now the Chow ring of a variety
pulls back isomorphically to the Chow ring of an affine bundle over it,
as we see by viewing an affine bundle as the difference of two
projective bundles. It follows that
 the Chow ring of $(V-S)/G$ agrees with the Chow
ring of $GL(N+n)/(GL(N)\times G)$ for this choice of $V$ and $S$.
\end{remark}

\begin{remark}
\label{eg}
 Edidin and Graham generalized the above definition of $CH^*BG$,
thought of as the $G$-equivariant Chow ring of a point. Namely, they defined
the $G$-equivariant Chow ring of a smooth $G$-variety $X$ by
$$CH^i_GX=CH^i(X\times (V-S))/G$$
for any $(V,S)$ as above such that the codimension of $S$ in $V$
is greater than $i$ \cite{EG}.
\end{remark}

\begin{proof} (Theorem \ref{indep})
We have to show
 that in codimension $i<s$, the Chow
groups of $(V-S)/G$ are independent of the subset $S$ and the representation
$V$. The independence of $S$ just uses the basic exact sequence for
Chow groups:
$$CH_*Y \arrow CH_*X\arrow CH_*U\arrow 0.$$
Here $X$ is a variety, $Y$ a closed subset, and $U=X-Y$. Apply this
to $X=(V-S)/G$ and $U=(V-S')/G$,
where $S'$ is some larger $G$-invariant subset of codimension
$\geq s$; since $CH_*(S'-S)/G$ vanishes in dimensions greater than
$\mbox{dim }(V-S)/G-s$, the Chow groups of $(V-S)/G$ map isomorphically
to the Chow
groups of $(V-S')/G$ in codimension less than $s$.
So the Chow ring
of $(V-S)/G$ is independent of $S$ in the range we consider.

The independence of $V$ follows from the double fibration construction,
used for example by Bogomolov \cite{Bogomolov}. That is, consider any
two representations $V$ and $W$ of $G$ such that $G$ acts freely
outside subsets $S_V$ and $S_W$ of codimension $\geq s$ and such that
the quotients $(V-S_V)/G$ and $(W-S_W)/G$ exist as varieties. Then consider
the direct sum $V\oplus W$. The quotient variety
$((V-S_V)\times W)/G$ exists, being a vector bundle over $(V-S_V)/G$,
and likewise the quotient variety $(V\times (W-S_W))/G$ exists,
as a vector bundle over $(W-S_W)/G$. Independence of $S$ 
(applied to the representation $V\oplus W$) shows
that these two vector bundles have the same Chow ring  in degrees
less than $s$. Since the Chow ring of a variety $X$ pulls back isomorphically
to the Chow ring of any vector bundle over $X$ \cite{Fulton},
 $(V-S_V)/G$ and $(W-S_W)/G$ have the  same
Chow ring in degrees less than $s$.
\end{proof}

\begin{proof} (Theorem \ref{twodefs})
Let $A^*BG$ (for the duration
of this proof) denote the ring
of characteristic classes for principal $G$-bundles over quasi-projective
varieties, as defined in the statement of Theorem \ref{twodefs}.
We know that the Chow ring of quotient varieties $(V-S)/G$, when
they exist, is independent of $V$ and $S$ in degrees
less than the codimension of $S$. We have a natural homomorphism
$$A^*BG\arrow CH^*(V-S)/G$$ when the quotient variety $(V-S)/G$
is quasi-projective.
Since there are some
quotients $(V-S)/G$ with the codimension of $S$ arbitrarily large
which are quasi-projective varieties (Remark \ref{quot}, above), we
get a natural homomorphism
$$A^*BG\arrow CH^*BG,$$
by the definition of the latter groups. We will show that this
homomorphism is an isomorphism.

\begin{lemma}
Let $G$ be a reductive group, and let $s$ be a positive integer.
For any quasi-projective variety $X$ with a principal $G$-bundle $E$,
there is an affine-space bundle $\pi:X'\arrow X$ and a map $f:X'\arrow (V-S)/G$
such that $\pi ^*E\cong f^*F$, where $F$ is the obvious principal 
$G$-bundle over $(V-S)/G$. Here $V$ is
a representation $V$ of $G$ and $S$ is a closed subset of codimension $\geq s$
such that $G$ acts freely on $V-S$ and the quotient variety $(V-S)/G$
exists as a quasi-projective variety.
\end{lemma}

\begin{proof}
By Jouanolou's trick \cite{Jouanolou}, for any quasi-projective
variety $X$ there is an affine-space bundle $\pi :X'\arrow X$ such that
$X'$ is affine. Pick any such $X'$ associated to the given variety $X$.
Then $X'$ is an affine variety with a principal $G$-bundle $\pi ^*E$; that is,
we have $X'=Y'/G$ for an affine variety $Y'$ with free $G$-action.

Let $Y'\subset A^n$ be any embedding of $Y'$ in affine space, given
by a finite-dimensional linear space $W_0\subset O(Y')$. Let $W\subset O(Y')$
be the $G$-linear span of $W_0$; then $W$ is a finite-dimensional
(by \cite{Mumford}, pp.~25-26)
representation of $G$ which gives a $G$-equivariant embedding 
$Y'\inj W^*$. Since $G$ acts freely on $Y'$, the orbits of $G$ on points
of $Y'$ are closed in $Y'$ and hence in $W^*$. That is, $Y'$ lies
in the $G$-stable subset of $W^*$ in the sense of geometric invariant
theory, which applies since $G$ is reductive
 \cite{Mumford}. The theory shows that 
the $G$-stable set is
an open subset of $W^*$ on which $G$ acts properly, so there is
 a smaller open subset $W^*-S\subset W^*$  containing $Y'$ on which
$G$ acts freely. Moreover,
the geometric quotient
$(W^*-S)/G$ exists as a quasi-projective variety \cite{Mumford}.
By adding another
representation of $G$ to $W^*$ if necessary, we can arrange that
$S$ has codimension at least any given number $s$. Thus we have
a map $Y=Y'/G\arrow (V-S)/G$ with the desired properties.
\end{proof}

It follows immediately from the lemma that the homomorphism
$A^*BG\arrow CH^*BG$
is injective. Namely, if $\alpha$ is a
characteristic class in $A^iBG$ which maps to 0 in $CH^iBG$, then
by definition of the latter group, $\alpha$
gives 0 on the quasi-projective quotients
$(V-S)/G$ with $S$ of codimension greater than $i$. 
By pulling back, the characteristic class gives 0 on
all smooth 
affine varieties with principal $G$-bundles. Since the Chow ring
of a smooth variety maps isomorphically to that of an affine bundle over it,
the characteristic class gives 0 on arbitrary  smooth quasi-projective
varieties with principal $G$-bundles. 
 By our definition of $A^*BG$, this means that $\alpha=0$.

To prove that $A^iBG\arrow CH^iBG$ is surjective, we argue as follows.
Given any element $\alpha \in CH^iBG$, which gives us an element
of $CH^i(V-S)/G$ for all quotient varieties $(V-S)/G$ with
$S$ of codimension greater than $i$,
we want to produce an element of $CH^iX$ for arbitrary
smooth quasi-projective varieties $X$ with $G$-bundles $E$. The lemma explains
how to do this: choose an affine-space bundle $\pi :X'\arrow X$ such that
there is a $G$-equivariant embedding $f:Y'\arrow V$ (where $X'=Y'/G$)
for some representation
$V$ of $G$ which is free   outside a subset $S$ of codimension greater
than $i$,
and pull back $\alpha \in CH^i(V-S)/G$ to $CH^iX'\cong CH^iX$. The only 
problem is to show that this is well-defined. We first show that it is
well-defined for a given affine $Y'$, independent of the choice
of the embedding $f$. Indeed, given any two embeddings $f:Y'\arrow  V$
and $g:Y'\arrow W$, we can consider the diagonal embedding
$h:Y'\arrow V\oplus W$; the commutative diagram
$$\begin{CD}
CH^i((V\oplus W)-S_{V\oplus W})/G  @= CH^i(V-S_V)/G\\
@V\cong VV    @VVV\\
CH^i(W-S_W)/G  @>>> CH^iY'/G
\end{CD}$$
shows that the two embeddings give the same element $\alpha \in CH^i Y'/G$.
Then, at last, if we have two different affine bundles over
the same quasi-projective variety $X$ which are both affine varieties,
there is a product affine
bundle which is an affine bundle over each one of them and in particular
is an affine variety. It follows that $\alpha \in CH^iX$ is well-defined
for each smooth quasi-projective variety $X$ with a principal $G$-bundle.
That is, $A^*BG\arrow CH^*BG$ is surjective as well as injective.
\end{proof}

\section{The factorization $CH^*BG\arrow MU^*BG\otimes_{MU^*}\Z
\arrow H^*(BG,\Z)$}
\label{general}

By \cite{TotaroMU}, for every smooth complex
algebraic variety $X$, the cycle map
$CH^*X\arrow H^*(X,\Z)$ factors naturally through a more refined
topological invariant of $X$, the ring $MU^*X\otimes_{MU^*}\Z$.
I will briefly describe what this means before explaining how it
extends to classifying spaces $BG$.

Namely, complex cobordism $MU^*X$ is a cohomology theory in the topological
sense; it is a graded ring associated to any topological space $X$
(see Stong \cite{Stong} for a general reference). When $X$ is a complex
manifold of complex dimension $n$, there is a simple interpretation
of the group $MU^iX$ as the bordism group of real $(2n-i)$-manifolds $M$ with
a complex linear structure on the stable tangent bundle $TM\oplus \R^N$
and with a continuous proper map $M\arrow X$. Another important
way to think of complex
cobordism is as the universal complex-oriented cohomology theory, where
a cohomology theory $h^*$ is complex-oriented if we are given an element
$c_1$ in the reduced cohomology
 $\tilde{h}^*(\C\Proj^{\infty})$ whose restriction to $\tilde{h}^*
\C\Proj^1$ freely generates the $h^*$-module $\tilde{h}^* (\C\Proj^1)$.
The familiar cohomology theories, integral cohomology and topological
$K$-theory, are complex-oriented in a natural way, which gives natural
transformations $MU^*X\arrow H^*(X,\Z)$ and $MU^*X\arrow K^*X$.
When $X$ is a point, the ring $MU^*:=MU^*(\text{point})$
 is a polynomial ring over $\Z$
on infinitely many generators, one in each degree $-2i$ for $i\geq 1$
(corresponding to some smooth compact complex $i$-fold, in terms
of the bordism definition of $MU^*$).
The map $MU^*\arrow H^*(\text{point},\Z)=\Z$ is the identity in degree 0
and (of course) 0 in degrees less than 0, whereas the map
$MU^*\arrow K^*$ takes the class of a complex $i$-fold in $MU^{-2i}$
to its Todd genus in $K^{-2i}=\Z$. We will need the map to $K$-theory
in section \ref{K_0}.

Using the trivial homomorphism $MU^*\arrow \Z$ coming from integral
cohomology, we can form the tensor product ring $MU^*X\otimes_{MU^*}\Z$;
this is a quotient ring of $MU^*X$ which maps to integral cohomology
$H^*(X,\Z)$. The map is an isomorphism for all spaces $X$ with
torsion-free integral cohomology, but not in general (as we will see
for many classifying spaces $BG$),
although it always becomes an isomorphism after tensoring
with the rationals. The main result of \cite{TotaroMU} is
that, for any smooth complex algebraic variety $X$, the cycle map
$CH^*X\arrow H^*(X,\Z)$ factors naturally through the ring
$MU^*X\otimes_{MU^*}\Z$. The factorization maps each codimension $i$
subvariety $Z$ in $X$
to the class in $MU^{2i}X$ of any resolution of singularities
for $Z$, using the above bordism interpretation of $MU^{2i}X$. This
class depends on the choice of resolution, but its image in
$MU^*X\otimes_{MU^*}\Z$ does not, and is invariant under rational
(or algebraic) equivalence of cycles.

Applying this to the varieties $(V-S)/G$ for representations $V$ of a complex
algebraic group $G$, we get a natural map $CH^*BG \arrow \invlim MU^*(V-S)/G
\otimes_{MU^*}\Z$, the inverse limit being taken over the 
representations $V$ of $G$. One can identify this inverse limit
with $MU^*BG\otimes_{MU^*}\Z$, giving the following theorem. Unfortunately,
this identification of the inverse limit is quite technical to prove,
and so this section might well be skipped on first reading.

\begin{theorem}
\label{factor}
For any complex algebraic group $G$,
the ring homomorphism $CH^*BG\arrow H^*(BG,\Z )$ factors naturally through
$MU^*BG\otimes_{MU^*}\Z$.
\end{theorem}

\begin{proof}
Here
 $MU^*BG$ must be viewed as a topological abelian group using
Landweber's  description of it
as $MU^*BG=\invlim MU^*(BG)_n$, where $(BG)_n$ denotes the
$n$-skeleton of a CW-complex representing $BG$ (see Lemma \ref{land},
below).  The tensor product
$MU^*BG\otimes_{MU^*}\Z$ must be defined as a topological tensor product:
$$MU^*BG\otimes_{MU^*}\Z:= \invlim (\im(MU^*BG \arrow MU^*(BG)_n)\otimes_{MU^*}\Z
),$$
as in \cite{TotaroMU}, p.~481. And what is needed to prove the lemma
is to show that, for any compact Lie group $G$, the ring
 $MU^*BG\otimes_{MU^*}\Z$ maps isomorphically
to the inverse limit
$$\invlim (MU^*(BG)_n\otimes_{MU^*}\Z).$$
Indeed, this last inverse limit is clearly the same as the inverse limit
of the rings $MU^*(V-S)/G\otimes_{MU^*}\Z$, which is the ring to which
$CH^*BG$ evidently maps, as explained above. (It is convenient to
casually identify a compact Lie group with its complexification here;
this is justified since the two groups have homotopy equivalent
classifying spaces.)

We have only a rather elaborate proof that $MU^*BG\otimes_{MU^*}\Z$
maps isomorphically to the above inverse limit. In fact, it seems that one can
show that the ring $MU^*BG\otimes_{MU^*}\Z$ is a finitely generated
abelian group in each degree, which means that these inverse systems
are in some sense trivial. We will not need that result in this paper.
More strongly, we can hope that $MU^*BG\otimes_{MU^*}\Z$
is a finitely generated $\Z$-algebra with generators in positive degree.

We need Landweber's general results on the complex cobordism of classifying
spaces  \cite{Landweber}:

\begin{lemma}
\label{land}
Let  $G$ be a compact Lie group, and let $(BG)_n$ denote
the $n$-skeleton of a CW complex $BG$.
Then
$$MU^*BG=\invlim MU^*(BG)_n.$$
Moreover, for each positive integer $n$, there is an $N\geq n$
such that, in all degrees at once,
$$\im (MU^*BG\arrow MU^*(BG)_n)=\im (MU^*(BG)_N\arrow MU^*(BG)_n).$$
\end{lemma}

To prove Theorem \ref{factor}, we need to show that the natural map
$$MU^*BG\otimes_{MU^*}\Z\arrow \invlim (MU^*(BG)_n\otimes
_{MU^*}\Z)$$
is an isomorphism.

Injectivity is fairly easy. Suppose that $x\in MU^*BG\otimes_{MU^*}\Z$ maps to
0 in the above inverse limit. By definition, $x$ itself is a compatible
sequence of elements $x_n\in \im(MU^*BG\arrow MU^*(BG)_n)\otimes_{MU^*}\Z$.
Lift $x_n$ to an element of the same name in $\im(MU^*BG\arrow MU^*(BG)_n)$.
We are assuming that, for each $n$, $x_n$ belongs to $MU^{<0}\cdot
MU^*(BG)_n$, and we want to show that for each $n$, $x_n$ belongs
to $MU^{<0}\cdot \im(MU^*BG\arrow MU^*(BG)_n)$.

By Lemma \ref{land}, for each $n$ there is an $N\geq n$ such that
$$\im(MU^*BG\arrow MU^*(BG)_n)=\im(MU^*(BG)_N\arrow
MU^*(BG)_n),$$
in all degrees at once.
 By hypothesis, we  know that $x_N$ belongs to $MU^{<0}\cdot MU^*(BG)_N$,
and we can restrict this equality to $(BG)_n$. We deduce the
desired conclusion about $x_n$, using the compatibility between
$x_N$ and $x_n$. So the above map is injective.

Now we only  need  to prove surjectivity of the homomorphism
$$MU^*BG\otimes_{MU^*}\Z\arrow \invlim (MU^*(BG)_n\otimes
_{MU^*}\Z).$$
Suppose we are given an element of the inverse limit on the right. Let
$x_n\in MU^*(BG)_n$, $n\geq 1$, be a choice of element representing
the given element of $MU^*(BG)_n\otimes_{MU^*}\Z$.
By Lemma \ref{land} again, for each $n$ there is an $N(n)\geq n$ such that
$$\im(MU^*BG\arrow MU^*(BG)_n)=\im(MU^*(BG)_{N(n)}\arrow
MU^*(BG)_n),$$
in all degrees at once.
Clearly we can assume that $N(1)< N(2)<\cdots $.
Let $y_n$ be the restriction  of $x_{N(n)}$ to $MU^*(BG)_n$;
by definition of $N(n)$, $y_n$ lifts to an element of
$MU^*BG$, which we also call $y_n$.

Moreover, we know a certain compatibility among the elements $y_n
\in MU^*BG$: all the elements $y_m$ for $m\geq n$ have the same restriction
to the ring $MU^*(BG)_{n}\otimes _{MU^*}\Z$, since
the $x$'s form an inverse system in these tensor product rings.
Equivalently, for any positive integers $m\geq n$, we can write
the restriction of $y_m-y_n$ to $MU^*(BG)_n$ as a finite sum
of products of elements of $MU^{<0}$ with elements of $MU^*(BG)_n$.
For each positive integer $n$, apply this fact to the positive integers
$N(n+1)\geq N(n)$: we find that the restriction of $y_{N(n+1)}-y_{N(n)}$
to $MU^*(BG)_{N(n)}$ belongs to $MU^{<0}\cdot MU^*(BG)_{N(n)}$,
meaning by this the additive group generated by such products. Then
the definition of $N(n)$ implies that upon further restricting
to $MU^*(BG)_n$, the element $y_{N(n+1)}-y_{N(n)}$
belongs to $MU^{<0}\cdot \im(MU^*BG\arrow MU^*(BG)_n)$.

Therefore, we can define one last sequence of elements  of
$MU^*BG$, $z_1,z_2,\dotsc$, by $z_1=y_{N(1)}$, $z_2$ equals $y_{N(2)}$
minus some element of $MU^{<0}\cdot MU^*BG$ 
such that $z_1$ and $z_2$
have the same restriction to $(BG)_1$, $z_3$ equals $y_{N(3)}$
minus some element of $MU^{<0}\cdot MU^*BG$ 
such that $z_3$ and $z_2$
have the same restriction to $(BG)_2$, and so on. Thus the $z_i$'s
form an element of $\invlim MU^*(BG)_n=MU^*BG$ (the equality
by Lemma \ref{land}). The resulting element maps to the element
of $\invlim (MU^*(BG)_n\otimes_{MU^*}\Z)$ we started with.
Theorem \ref{factor} is proved.
\end{proof}

\section{$K_0BG$, and the Chow ring with small primes inverted}
\label{K_0}

One can define a ring $K_0BG$ analogous to $CH^*BG$ using
the algebraic $K$-group $K_0$ in place of the Chow ring.
The ring $K_0BG$ can
be completely computed using a theorem of Merkurjev \cite{Merkurjev}.
 The result is
an algebraic analogue of the theorem of Atiyah-Hirzebruch-Segal on
the topological $K$-theory of classifying spaces \cite{Atiyah}, \cite{AS}.
We will use the calculation of $K_0BG$ to derive some good information
on $CH^*BG$ and $MU^*BG\otimes_{MU^*}\Z$ with small primes inverted. 

For any linear algebraic group $G$ over a field $k$, we define
$$K_0BG=\invlim K_0(V-S)/G,$$
where the inverse limit runs over representations $V$ of $G$ over $k$.
Here $K_0X$ is the Grothendieck
group of algebraic vector bundles over a variety. This (standard)
 notation has the
peculiarity that for a variety $X$ over the complex numbers, there is
natural map from the algebraic $K$-group $K_0X$ to the topological $K$-group
$K^0X$.

To state the calculation of $K_0BG$, for any algebraic group $G$
over a field $k$, let $R(G)$ denote the representation ring of $G$ over $k$
(the Grothendieck group of finite-dimensional representations of $G$ over $k$).
Thus $R(G)$ is the free abelian group on the set of irreducible representations
of $G$ over $k$. Let $R(G)\sphat$ be the completion of $R(G)$ with respect
to powers of the augmentation ideal $\text{ker}(R(G)\arrow \Z)$.

\begin{theorem}
\label{merk}
For any algebraic group $G$ over a field $k$, there is an isomorphism
$$R(G)\sphat \arrow K_0BG.$$
\end{theorem}

For comparison, Atiyah, Hirzebruch, and Segal proved that for a complex
algebraic group $G$ (or, equivalently, a compact Lie group), the completed
representation ring $R(G)\sphat$ maps isomorphically to the topological
$K$-group $K^0BG$, whereas $K^1BG=0$ \cite{Atiyah}, \cite{AS}.
Thus Theorem \ref{merk} implies
that for a complex algebraic group $G$, the algebraic $K$-group $K_0BG$
maps isomorphically to the topological $K$-group $K^0BG$.

\begin{proof}
Merkurjev computed, in particular, the algebraic $K$-group $K_0$ for
 any homogeneous space
of the form $GL(n)/G$:
$$K_0(GL(n)/G)=\Z\otimes_{R(GL(n))}R(G)$$
\cite{Merkurjev}. By Remark \ref{quot}, applied to $K_0$ rather than the Chow
ring, we can compute $K_0BG$ as
$$K_0BG=\invlim_N K_0(GL(N+n)/GL(N)\times G)$$
for any faithful representation $G\arrow GL(n)$. So, by Merkurjev's
theorem, we have
\begin{align*}
K_0BG &= \invlim_N \Z\otimes_{R(GL(N+n))}R(GL(N)\times G)\\
 &= \invlim_N [\Z\otimes_{R(GL(N+n))}R(GL(N)\times GL(n))]\otimes_{R(GL(n))}
R(G).
\end{align*}
The ring in brackets is $K_0(GL(N+n)/GL(N)\times GL(n))$, or equivalently
$K_0$ of the Grassmannian of $n$-planes in affine $(N+n)$-space. In particular,
we easily compute that the representation ring of $GL(n)$ maps onto this ring,
with kernel contained in $F^{N+1}_{\gamma}R(GL(n))$. (See Atiyah \cite{Atiyah}
for the definition of the gamma filtration of a $\lambda$-ring.) It follows
that the representation ring $R(G)$ maps onto the above ring
$K_0(GL(N+n)/GL(N)\times G)$ with kernel contained in $F^{N+1}_{\gamma}R(G)$.
Also, this kernel contains $F^M_{\gamma}R(G)$ for some large $M$, just because
$GL(N+n)/GL(N)\times G$ is finite-dimensional. It follows that
$K_0BG$ is isomorphic to the completion of the representation ring
$R(G)$ with respect to the gamma filtration. By Atiyah's arguments
(\cite{Atiyah}, pp.~56-57), this is the same as the completion of $R(G)$
with respect to powers of the augmentation ideal.
\end{proof}

Esnault, Kahn, Levine, and Viehweg pointed out that Merkurjev's
theorem implies that $CH^1BG$ and $CH^2BG$ are generated by Chern classes
of representations of $G$ (\cite{EKLV}, Appendix C). In general,
we have the following statement.

\begin{corollary}
\label{lowcor}
For any algebraic group $G$ over a field $k$, the subgroup of the
Chow group $CH^iBG$ generated by Chern classes $c_i$ of $k$-representations
of $G$ contains $(i-1)!CH^iBG$.
\end{corollary}

\begin{proof}
For any smooth
algebraic variety $X$ over a field $k$,
 the algebraic $K$-group $K_0X$ has two natural
filtrations, the geometric filtration (by codimension of support)
and the gamma filtration \cite{Grothint}. There are natural maps
from both $CH^iX$ and $\gr_{\gamma}^iK_0X$ to $\gr_{\text{geom}}^iK_0X$.
The map $CH^iX\arrow \gr_{\text{geom}}^ iK_0X$ is surjective, and
the $i$th Chern class gives a map back such that the composition
$$\begin{CD}
CH^iX @>>> \gr_{\text{geom}}^iK_0X @>c_i>> CH^iX
\end{CD}$$
is multiplication by $(-1)^{i-1}(i-1)!$, by Riemann-Roch
without denominators \cite{Jouanolou}, \cite{Grothint}. It follows
that the surjection from $CH^iX$ to $\gr_{\text{geom}}^iK^ 0X$
becomes an isomorphism after inverting $(i-1)!$, and the subgroup of
the Chow group $CH^iX$ generated by Chern classes $c_i$ of
vector bundles on $X$ contains $(i-1)!CH^iX$.

Apply this to the variety
$X=GL(N+n)/(GL(N)\times G)$, for any faithful representation $G\arrow GL(n)$
and any $N\geq i$. By Remark \ref{quot}, $CH^iBG$ maps isomorphically
to $CH^iX$, and by Merkurjev's theorem, $K_0X$ is generated by representations
of $G$. It follows that the subgroup of $CH^iBG$ generated by Chern classes
$c_i$ of representations of $G$ contains $(i-1)!CH^iBG$.
\end{proof}

Now suppose that $G$ is an algebraic group over the complex numbers.
Then, in keeping with the observation that the homomorphism
$CH^*BG\arrow MU^*BG\otimes_{MU^*}\Z$ is often an isomorphism,
we can prove properties analogous to the above property of $CH^*BG$
for $MU^*BG\otimes_{MU^*}\Z$, as follows.

\begin{theorem}
\label{smalltors}
For any complex algebraic group $G$, the subgroup of $CH^iBG$ generated
by Chern classes $c_i$ of representations of $G$ contains $(i-1)!CH^iBG$,
and the subgroup of $(MU^*BG\otimes_{MU^*}\Z)^{2i}$ generated by Chern
classes $c_i$ of representations contains $(i-1)!(MU^*BG\otimes_{MU^*}\Z)
^{2i}$.
 
Moreover, there are maps
$$\begin{CD}
\gr^i_{\gamma}R(G)@>c_i>> CH^iBG @>>> (MU^*BG\otimes_{MU^*}\Z)
^{2i} @>>> \gr^{2i}_{\text{top}}R(G)
\end{CD}$$
such that the first two maps are surjective after inverting $(i-1)!$.
(We will give Atiyah's definition of the topological filtration on $R(G)$
in the proof.)
The composition of all three maps is $(-1)^{i-1}(i-1)!$ times
the natural map from $\gr^i_{\gamma}R(G)$ to $\gr^{2i}_{\text{top}}
R(G)$. So, for example, if that natural map is injective, 
then the groups $\gr^i_{\gamma}R(G)$, $CH^iBG$, and
$(MU^*BG\otimes_{MU^*} \Z)^{2i}$ all become isomorphic after inverting
$(i-1)!$.

Finally, the odd-degree group $(MU ^*BG\otimes_{MU^*}\Z)^{2i+1}$
is killed by $i!$.
\end{theorem}
 
\begin{proof}
For any smooth algebraic variety $X$, the proof of Corollary \ref{lowcor}
shows that the natural map
$$CH^iBG\arrow \gr^i_{\text{geom}}K_0X$$
and the Chern class map
$$c_i: \gr^i_{\text{geom}}K_0X \arrow CH^iBG$$
become isomorphisms after inverting $(i-1)!$. In particular, both groups
are generated by Chern classes of elements of $K_0X$ after inverting $(i-1)!$.
Also, by definition of the gamma filtration,the image of
the natural map
$$\gr^i_{\gamma}K_0X\arrow \gr^i_{\text{geom}}K_0X$$
is the subgroup of $\gr^i_{\text{geom}}K_0X$ generated by Chern classes
$c_i$ of elements of $K_0X$. It follows that this map becomes surjective,
although not necessarily injective, after inverting $(i-1)!$.
So the composition
$$\begin{CD}
\gr^i_{\gamma}K_0X@>>> \gr^i_{\text{geom}}K_0X
@>c_i>> CH^iX
\end{CD}$$
becomes surjective after inverting $(i-1)!$.

We can apply all this to the  classifying space $BG$, viewed as  a limit
of smooth algebraic varieties over $k$,   using the identification
of the algebraic $K$-group $K_0BG$ with the completed representation
ring of $G$ (Theorem \ref{merk}). The map $\gr^i_{\gamma}R(G)
\arrow CH^iBG$ in the statement is defined as the composition
$$\begin{CD}
\gr^i_{\gamma}R(G)@>>> \gr^i_{\text{geom}}R(G)
@>c_i>> CH^iBG.
\end{CD}$$
It becomes surjective after inverting $(i-1)!$, as we want.

We showed in Corollary \ref{lowcor} that the subgroup of $CH^iBG$ generated
by Chern classes $c_i$ of representations of $G$ contains $(i-1)!CH^iBG$.
Now let us prove the analogous statement for
complex cobordism. We need the following ``Riemann-Roch without
denominators''
formula for complex cobordism. The statement uses the map of cohomology
theories $MU^*X\arrow K^*X$ defined by Atiyah-Hirzebruch and Conner-Floyd
\cite{AHriem}, \cite{CF}; when $X$ is a point, the homomorphism $MU^*
\arrow K^*$ is simply the Todd genus.
Using the periodicity of topological $K$-theory, we have
a map $MU^{2i}X\arrow K^{2i}X=K^0X$.

 \begin{lemma}
\label{rr}
For any topological space $X$, the composition
$$\begin{CD}
MU^{2i}X @>>> K^0X @>c_i>> (MU^*X\otimes_{MU^*}\Z)^{2i}
\end{CD}$$
is multiplication by $(-1)^{i-1}(i-1)!$.
\end{lemma}

\begin{proof}
Let us first show that this is correct after we map further from the group
$(MU^*X\otimes_{MU^*}\Z)^{2i}$ into rational cohomology $H^{2i}(X,\Q)$.
It will suffice to show that the image in $K^0X$ of an element $x\in MU^{2i}X$
has Chern character $\ch_i$ in $H^{2i}(X,\Q)$ equal to the image of $x$ under
the obvious map and has $\ch_j=0$ for $j<i$. This is enough because
the Chern character $\ch_i$ of an element of $K^0X$ for which
$ch_j$ is 0 for $j<i$ is given by
$$\text{ch}_i=(-1)^{i-1}c_i/(i-1)!.$$

Because complex cobordism is a connective cohomology theory
($MU^*(\text{point})$ is 0 in positive degrees), any element $x\in MU^{2i}X$
restricts to 0 on the $2j$-skeleton of X for $j<i$, which implies
that the image of $x$ in $K^0X$ has Chern character $\ch_j$ equal to 0
in $H^{2j}(X,\Q)$ for $j<i$. Also, because $MU^0(\text{point})$ is just $\Z$,
the natural transformation from $MU^{2i}X$ to $H^{2i}(X,\Q)$ given by mapping
to $K^0X$ and taking the Chern character $\ch_i$
must be a rational
multiple of the obvious map. We want to show that the rational multiple
is 1. It suffices to check this when $X$ is the pair $(S^{2i},\text{point})$.
Then $MU^{2i}X=MU^0=\Z$ via the suspension isomorphism,
and since the map from complex cobordism to $K$-theory is multiplicative,
the element 1 in this group maps to 1 in $K^{2i}X$ when we identify
this group with $K^0=\Z$ via the suspension isomorphism. So the calculation
comes down to: if we start with the element $1$ in $K^0=\Z$, identify it
with an element $K^{2i}X$ (where $X=(S^{2i},\text{point})$)
via the suspension isomorphism, and then
identify this with $K^0X$ via periodicity, show that the Chern character
$\ch_i$ of the resulting element is 1 in $H^{2i}(X,\Q)=\Q$. This is
a classical property of Bott periodicity (\cite{AHK}, p.~16).

To prove the above formula as stated, it suffices to prove it
in the universal case, where $X$ is the $2i$th space in
the $\Omega$-spectrum $MU$.  The point is that $X$ has torsion-free
cohomology by Wilson (\cite{Wilson1}, pp.~52-53). So $MU^*X\otimes_{MU^*}\Z$
is equal to $H^*(X,Z)$ (by the Atiyah-Hirzebruch spectral sequence,
as in \cite{TotaroMU}, p.~471) and injects into $H^*(X,\Q)$. Thus
the desired formula is true since it is true in rational cohomology.
\end{proof}

Applying Lemma \ref{rr} to the classifying space of a complex algebraic
group $G$, we deduce that the subgroup of $MU^*BG\otimes _{MU^*}\Z$ in degree
$2i$ which is generated by  Chern classes $c_i$ of representations  of $G$
contains $(i-1)!$ times the whole group. Since
the topological $K$-group $K^1BG$ is $0$,
applying Lemma \ref{rr} to the suspension of $BG$ shows that
$MU^*BG\otimes_{MU^*}\Z$ in odd  degree $2i+1$ is killed by $i!$.

Finally, for any CW complex $X$, the map from $MU^{2i}X$ to $K^0X$
lands in $F^{2i}_{\text{top}}K^0X$, where the topological filtration
of $K^0X$ is defined by letting $F^{2i}_{\text{top}}K^0X$ be the subgroup
of $K^0X$ which restricts to 0 on the $(2i-1)$-skeleton of $X$.
This is clear because every element of $MU^{2i}X$ restricts to 0 on
the $(2i-1)$-skeleton (because complex cobordism, unlike $K$-theory,
is a connective cohomology theory, meaning that $MU^*(\text{point})$ is
0 in positive degrees).
The quotient group $( MU^*X\otimes_{MU^*}\Z)^{2i}$ maps to
the associated graded group $\gr^{2i}_{\text{top}}K^0X$, and it is
elementary to identify the composition
$$\begin{CD}
\gr^i_{\gamma}K^0X@ >c_i>> (MU^*X\otimes_{MU^*}\Z)^{2i} @>>>
\gr^{2i}_{\text{top}}K^0X
\end{CD}$$ 
with $(-1)^{i-1}(i-1)!$ times the natural map. For $X=BG$, this completes
the proof of Theorem \ref{smalltors}, since the topological filtration
on the representation ring $R(G)$ is defined as the pullback of
the topological filtration of $K^0BG$.
\end{proof}

\begin{corollary}
\label{lowdeg}
For any complex algebraic group $G$,
the map from the Chow group $CH^iBG$ to $(MU^*BG\otimes_{MU^*}\Z)^{2i}$ is an
isomorphism for $i\leq 2$. The map from these groups to
$H^{2i}(BG,\Z)$ is an isomorphism for $i\leq 1$ and injective
for $i=2$. Also, for $i\leq 2$, these groups are generated by
Chern classes $c_i$ of representations of $G$.
\end{corollary}

\begin{proof}
The map from $CH^iBG$ to $(MU^*BG\otimes_{MU^*}\Z)^{2i}$
 is an isomorphism for $i=1$ by Theorem \ref{smalltors}, using
that $\gr^1_{\gamma}R(G)=\gr^1_{\text{top}}R(G)=H^2(BG,\Z)$
(\cite{Atiyah}, p.~58). Unfortunately, this argument is not enough
to prove the desired isomorphism for $i=2$, since E.~Weiss
found for the group $G=A_4$ that the map  from $ \gr^2_{\gamma}R(G)$
to $\gr^2_{\text{top}}R(G)$ has a kernel of order 2 (see Thomas
\cite{Thomas}). At least Theorem \ref{smalltors} shows that
the map from  $CH^2BG$ to $(MU^*BG\otimes_{MU^*}\Z)^4$   is surjective
with finite kernel.

The point is to apply an argument of Bloch and Merkurjev-Suslin,
as formulated by Colliot-Th\'el\`ene (\cite{Colliot}, p.~13). For each
smooth algebraic variety $X$ over a field $k$ and each prime number
$l$ invertible in $k$, this argument gives
a commutative diagram
$$\begin{CD}
H^1_{\text{Zar}}(X,\mathcal{H}^2_{\text{et}}(\Q_l/\Z_l(i))) @>>\gamma_2 >
H^3_{\text{et}}(X,\Q_l/\Z_l(i))\\
@V\alpha_2VV @VVV \\
CH^2(X)_{l{\text{-tors}}} @>>> H^4_{\text{et}}(X,\Z_l(i))
\end{CD}$$
such that $\gamma_2$ is injective and $\alpha_2$ is surjective.
Take $X$ to be a quotient variety $(V-S)/G$ for a complex algebraic
group $G$, with $S$ of large codimension in $V$. Then $H^3(X,\Q_l/\Z_l)$
maps injectively to $H^4(X,\Z_l)$, since $H^3(X,\Q)=H^3(BG,\Q)=0$.
It follows that $CH^2(X)_{l{\text{-tors}}}$ injects into
$H^4(X,\Z_l)$, or equivalently that $CH^2(BG)_{l{\text{-tors}}}$
injects into $H^4(BG,\Z_l)$. Since we know that the map
from $CH^2BG$ to $(MU^*BG\otimes_{MU^*}\Z)^4$ is surjective with
finite kernel, it must be an isomorphism, and these groups
must inject into $H^4(BG,\Z)$.

Finally, Corollary \ref{lowcor} shows that $CH^2BG$ is generated by
Chern classes $c_2$ of representations. 
\end{proof}

It follows,
for example, that $CH^2BG\arrow H^4(BG,\Z)=\Z$ is not surjective,
for all simply connected simple groups $G$ over $\C$ other than $SL(n)$
and $Sp(2n)$. The analogous statement for $G=SO(4)$ was the
source of the examples of varieties with nonzero Griffiths group
in \cite{TotaroMU}.

\section{Transferred Euler classes}
\label{surjcases}

As mentioned in the introduction, for any complex reductive group $G$,
there is a natural class of elements
of     the ring $MU^*BG\otimes_{MU^*}\Z$ which lie in the image 
of the Chow ring of $BG$.
  For any $i$-dimensional   
representation $E$ of a subgroup $H\subset G$ (including $G$ itself),
 $E$ determines an algebraic vector bundle on $BH$ (meaning
an algebraic vector bundle on all the varieties  $(V-S)/H$),   so it has 
Chern classes in the Chow ring, in complex cobordism, and in
ordinary cohomology. These are compatible under the natural
homomorphisms
$$CH^*BH\arrow MU^*BH\otimes_{MU^*}\Z\arrow H^*(BH,\Z).$$
In particular, we define the Euler class $\chi(E)$ to mean the
top Chern class of $E$ in any of these groups.  These homomorphisms
are  also compatible with transfer maps, for $H$ of finite index in $G$,
 so the transferred
Euler class $\text{tr}_H^G \chi(E)$ is an element of $CH^*BG$
which maps to the element with the same name in $MU^*BG\otimes_{MU^*}\Z$
and in $H^*(BG,\Z)$.
Sometimes we also call any $\Z $-linear combination of transferred
Euler classes a transferred Euler class. In this terminology, if $G$
is a finite group, the
transferred Euler classes form a subring of any of the above
three rings,  by the
arguments of Hopkins, Kuhn, and Ravenel \cite{HKR}.

There are many finite groups $G$
which satisfy the conjecture
of Hopkins, Kuhn, and Ravenel that the Morava $K$-theories of $BG$
are generated by transferred Euler classes; by
Ravenel, Wilson, and Yagita, it then follows that
$MU^*BG$ is generated as a topological $MU^*$-module by transferred
Euler classes \cite{RWY}. Whenever $G$ has this property, the Chow ring of $BG$
clearly maps onto $MU^*BG\otimes_{MU^*}\Z$. For such groups $G$,
we have a computation of $\im (CH^*BG\arrow H^*(BG, \Z ))$
in topological terms: it is the $\Z $-submodule of $H^*(BG, \Z )$
generated by transferred Euler classes, or equivalently it is the image
of $MU^*BG\arrow H^*(BG, \Z )$. The conjecture has been checked for
various finite groups, including abelian groups, groups of order $p^3$,
the symmetric groups, and finite groups of Lie type
away from the defining characteristic, by \cite{TY}, \cite{HKR},
and \cite{Tanabe}. The conjecture about Morava $K$-theory fails
for the group $G$
considered by Kriz \cite{Kriz}. We can still guess that the Chow ring of $BG$
is additively generated by transferred Euler classes for all finite
groups $G$.

As the reader has probably noticed,
we could consider the
transfers of arbitrary Chern classes of representations
instead of just the top Chern class.
But one can show that
transfers of arbitrary Chern classes are $\Z $-linear combinations of
transfers of Euler classes, so that this would not give anything new.
On the other hand we definitely cannot avoid mentioning transfers; that is,
$MU^*BG$ and $CH^*BG$ are not generated as algebras by the Chern classes
of representations of $G$ itself, for some groups $G$. Chris Stretch
showed in 1987 that $MU^*(BS _6)$ is not generated by Chern classes;
we will now give his pleasant argument.

The point is that $S _6$ has two subgroups which are pointwise
conjugate, but not conjugate:
$$H_1=\{1, (12)(34), (12)(56), (34)(56)\}$$
$$H_2=\{1, (12)(34), (13)(24), (14)(23)\}$$
Since they are pointwise conjugate, character theory shows that any
complex
representation of $S _6$ has a restriction to $H_1$ which is isomorphic
to its restriction to $H_2$
(in terms of any fixed isomorphism $H_1\cong H_2$). It follows that any element
of $MU^*(BS _6)$ in the subring generated by Chern classes of
representations has the same restriction to $H_1$ as to $H_2$, in terms of
the fixed identification $H_1\cong H_2$. So it suffices to find an element
of $MU^*(BS _6)$ which restricts differently to the two subgroups.

One such element is the transfer from $S_4\times S_2\subset S_6$
to $S_6$ of the Euler class of the representation
$S_4\times S_2\arrow S_4\arrow GL(3,\C)$. (The last representation is
the permutation representation of $S_4$ minus the trivial representation.)
The double coset formula gives the restriction of this transferred element
to any subgroup, and one finds that the restrictions of this element
of $MU^*BS _6$ to $MU^*BH_1$ and $MU^*BH_2$ are different, since their
 images in $H^*(BH_1,\Z)$ and $H^*(BH_2,\Z)$
are different. (This argument actually
gives the stronger conclusion that the image of $MU^*BS_6\arrow
H^*(BS_6,\Z)$ is not generated by Chern classes.)

\section{Chow groups of Godeaux-Serre varieties}
\label{godeaux}

For any finite group $G$, Godeaux and Serre constructed smooth
projective varieties with fundamental group $G$, as quotients
of complete intersections by free actions of $G$ (\cite{Serretop},
section 20). Atiyah and Hirzebruch showed the falsity of
the Hodge conjecture for integral cohomology using some
of these varieties \cite{Atiyah-Hirzebruch}: not all the torsion cohomology
classes on these varieties can be represented by algebraic cycles.
So it becomes a natural problem
to compute the Chow groups of these varieties. We can at least
offer a conjecture. Let $G$ be a finite group, $V$ a representation
of $G$ over a field $k$, and $X$ a smooth complete intersection
defined over $k$ in $P(V)$ such
that $G$ acts freely on $X$. We call the quotient variety $X/G$ a
{\it Godeaux-Serre variety.}

\begin{conjecture}
The natural homomorphism
$$CH^i(BG\times BG_m)_k\arrow CH^iX/G$$
is an isomorphism for $i<\text{dim }X/2$. For very general $X$
of sufficiently high degree, this homomorphism is an isomorphism
for $i<\text{dim }X$.
\end{conjecture}

The homomorphism mentioned in the conjecture comes from
the natural $G\times G_m$-bundle over $X/G$, where the $G$-bundle
 over $X/G$ is obvious and the $G_m$-bundle on $X/G$ is the one
which corresponds to the $G$-equivariant line bundle $O(1)$ on
$X\subset P(V)$.

The proof that Godeaux-Serre varieties exist \cite{Serretop} 
(for every finite group
$G$ and in every dimension $r$) helps to show why the conjecture is plausible.
Namely, let $G$ be a finite group and $V$ a representation  of $G$.
We can imbed the quotient variety $P(V)/G$ in some projective space
$\Proj^N$. Let $S'$ be the closed subset of $P(V)$ where $G$ does not
act freely; for any given $r\geq 0$, we can choose the representation
$V$ of $G$ such that $S'$ has codimension greater than $r$ in $P(V)$.
Then we construct Godeaux-Serre varieties $X/G$ of dimension $r$
by intersecting $P(V)/G$ with a general linear space in $\Proj^N$
of codimension equal to $\text{dim }P(V)/G -r$. The intersection will not
meet $S'/G$ and will be a smooth compact subvariety of the
smooth noncompact variety $(P(V)-S')/G$.
We can think of $X/G$ as a complete intersection in
$(P(V)-S')/G$.

Let $S\subset V$ be the union of 0 with the inverse image of
$S'\subset P(V)$; then $(P(V)-S')/G$ can also be viewed as the quotient
$(V-S)/(G\times G_m)$. 
Then $S$ has codimension in $V$ greater than $r=\text{dim }X$,
so $$CH^i(BG\times BG_m)=CH^i((V-S)/(G\times G_m))$$
for $i\leq r$. So the above conjecture would follow if we could
relate the Chow groups of $(V-S)/(G\times G_m)$ to those of the
complete intersection $X/G$ inside it. The precise bounds in
the conjecture are those suggested by the usual conjectures
on complete intersections in smooth projective varieties, due to
Hartshorne \cite{Hartshorne} and Nori (\cite{Nori}, p.~368);
see also Paranjape \cite{Paranjape}, pp.~643-644.
The analogous topological statement was proved in this case
by Atiyah and Hirzebruch: for Godeaux-Serre varieties
over the complex numbers, the natural homotopy class of maps
$X/G\arrow BG\times BS^1$ is  $r$-connected, where $r=\text{dim}_{\C}X$
(\cite{Atiyah-Hirzebruch}, p.~42).

For this conjecture to really say anything about the Chow groups
of Godeaux-Serre varieties, we need to understand the Chow ring
of $BG\times BG_m$.
We can think of $BG_m$ as the infinite-dimensional projective space,
which makes it easy to check that the Chow ring  of
$BG\times BG_m$ is a polynomial ring in one variable in degree 1 over the
Chow ring of $BG$. For general finite groups $G$, we have
no computation of the Chow ring of $BG$, but it tends to be computable,
as shown by various  results in this paper.

It remains to say what kind of evidence can be offered for the conjecture.
Even for the trivial group $G=1$, the conjecture is far
out of reach, since in that case it amounts to the conjectures
of Hartshorne and Nori for complete intersections in projective space.
For example, Hartshorne's conjecture says that the low-codimension
Chow groups of a smooth complete intersection over 
the complex numbers are equal to $\Z$,
but for all we know, they might even be uncountable (apart from $CH^1$).
Nonetheless, there are various suggestive pieces of evidence for
the above conjecture. Let us work over $\C$ in what follows.

A weak sort of evidence is that the only known examples of
Godeaux-Serre varieties for which $CH^*(X/G)\arrow H^*(X/G,\Z)$
can be proved to be non-surjective \cite{Atiyah-Hirzebruch} or non-injective
\cite{TotaroMU} in low codimension come from groups $G$ for which
the corresponding map for
$BG\times BG_m$ is non-surjective or non-injective.

Also, since we do not know how to prove the triviality of
the low-codimension Chow groups of a complete intersection $X$
in projective space, we might settle for trying to understand
the kernel of the pullback map $CH^*(X/G)\arrow CH^*X$. By
the obvious transfer argument, this kernel is a torsion group,
killed by $|G|$. The torsion subgroup of the Chow groups
of any variety over an algebraically closed field is countable
by Suslin (\cite{Suslin}, p.~227). In that sense, the kernel
of $CH^*(X/G)\arrow CH^*X$ is under better control than $CH^*X$.

For $CH^2$, we can prove the conjecture on this kernel for varieties $X/G$
of sufficiently large dimension:

\begin{theorem}
For any finite group $G$, let $X$ be a smooth complete intersection
over $\C$ on which $G$ acts freely. For $\text{dim }X$ sufficiently
large, depending on $G$, the natural map
$$\ker(CH^2(BG\times B\C^*)\arrow CH^2(B\C^*))\arrow
\ker(CH^2(X/G)\arrow CH^2X)$$
is an isomorphism. The first kernel is naturally identified with
$CH^1BG\oplus CH^2BG$.
\end{theorem}

\begin{proof}
By the argument  of Bloch and Merkurjev-Suslin used
already in section \ref{K_0},
the torsion subgroup of $CH^2$ of a smooth projective variety $Y$ over
$\C$ maps injectively to $H^3(Y,\Q/\Z)$ (see \cite{Colliot}, p.~17).
Taking $Y$ to be a Godeaux-Serre variety $X/G$ of dimension at least 4,
we have $H^3(X/G,\Q)\inj H^3(X,\Q)=0$, and so the torsion subgroup
of $CH^2(X/G)$ maps injectively
 to the ordinary cohomology group $H^4(X/G,\Z)$. If
$X$ has dimension at least 5, this cohomology group is equal to
$H^4(BG',\Z)$, where we write $G'=G\times \C^*$.
 Thus the kernel of $CH^2X/G\arrow CH^2X$  maps injectively to
$H^4(BG',\Z)$. The image in cohomology of that kernel clearly
contains the image of the kernel
of $CH^2BG'\arrow CH^2B\C^*=\Z$, and we will show that equality holds
when $X$ has sufficiently large dimension, depending on $G$.

In fact, we know by Corollary \ref{lowdeg} that the map 
$$CH^2BG'\arrow
(MU^*BG'\otimes_{MU^*}\Z)^{4}$$
  is an isomorphism. And if $(BG)_n$
denotes the $n$-skeleton of $BG$, Lemma \ref{land} (Landweber's
theorem) says that for every $n$
there exists an $r\geq n$ such that
$$\im(MU^*BG'\arrow MU^*(BG')_n)=\im(MU^*(BG')_r\arrow MU^*(BG')_n).$$
Now, for any space $X$, the image of an element $x\in MU^4X$ in $H^4(X,\Z)$
only depends on the restriction of $x$ to $MU^4$ of the 4-skeleton $X_4$,
since $H^4(X,\Z)$ injects into $H^4(X_4,\Z)$.
So the above statement implies that there is an $r\geq 4$ such that
$$\im(MU^4BG'\arrow H^4(BG',\Z))=\im(MU^4(BG')_r\arrow H^4(BG',\Z)).$$
Since a Godeaux-Serre variety $X/G$ of dimension $r$ contains the
$r$-skeleton of $BG'$ up to homotopy, it follows that
$$\im(MU^4BG'\arrow H^4(X/G,\Z))=\im(MU^4X/G\arrow H^4(X/G,\Z))$$
for Godeaux-Serre varieties $X/G$ of dimension at least $r$.
Thus, for Godeaux-Serre varieties of sufficiently large dimension,
the image of the Chow group $CH^2X/G$ in integral cohomology
is contained in the image of $MU^4X/G$, hence in the image of
$MU^4BG'$, and hence (by Corollary \ref{lowdeg}) in the image of $CH^2BG'$.
In particular, if we start with an element $x\in CH^2(X/G)$ which pulls back
to 0 in $CH^2X$, then, as we have said,
the class of $x$ in integral cohomology is the image
of an element of $CH^2BG'$, and this element must pull back to 0 in $CH^2B\C^*
=H^4(B\C^*,\Z)=\Z$ because $x$ pulls back to 0 in $H^4(X,\Z)=H^4(B\C^*,\Z)=\Z$.

Finally, the maps of both $\ker(CH^2BG'\arrow \Z)$
 and $\ker(CH^2X/G\arrow CH^2X)$
to $H^4(X/G,\Z)=H^4(BG',\Z)$ are injective. For the first group,
this is Corollary \ref{lowdeg}. For the second group, we proved this
two paragraphs back. It follows that the kernel of
$CH^2BG'\arrow CH^2B\C^*=\Z$
 maps isomorphically to the kernel of $CH^2X/G\arrow CH^2X$.
\end{proof}

\section{Some examples of the Chow K\"{u}nneth formula}
\label{Kunneth}

As a warmup to our discussion of the Chow groups of symmetric products
(which we will use to compute the Chow cohomology of the symmetric group),
we make some easy remarks about the Chow groups of a product in this section.

Let $X$ and $Y$ be algebraic varieties over a field $k$.
We describe a special situation in which the Chow groups of $X\times Y$
are determined by the Chow groups of $X$ and $Y$ by the ``Chow K\"{u}nneth
formula'': $$CH_*(X\times Y)\cong CH_*(X)\otimes _{\Z }CH_*Y.$$
This formula is certainly false in general, for example for the product
of an elliptic curve with itself.
The interesting thing is that the formula
holds in some cases even where $CH_*X$ and $CH_*Y$ have a lot of torsion,
whereas the corresponding formula for integer homology would have an
extra ``Tor'' term.

\begin{lemma}
If $X$ is any variety and $Y$ is  a variety which can be partitioned
into open subsets of affine spaces, then the map
$$CH_*X\otimes _{\Z }CH_*Y\arrow CH_*(X\times Y)$$
is surjective.
\end{lemma}

\begin{proof}
We recall the basic exact sequence for the Chow groups.
Let $Y$ be a variety, $S\subset Y$ a closed subvariety, and $U=Y-S$.
Then the sequence
$$CH_*S\arrow CH_*Y\arrow CH_*U\arrow 0$$
is exact.

We know the Chow groups of affine space: $CH_i(A^n)$ is $\Z $ if
$i=n$ and 0 otherwise. The basic exact sequence implies that
the Chow groups of a nonempty open subset $Y\subset A^n$ are the same:
$CH_*Y=\Z $ in dimension $n$. We also know that
$$CH_*X\otimes_{\Z }CH_*A^n\cong CH_*(X\times A^n)$$
for any variety $X$. It follows that
$$CH_*X\otimes _{\Z }CH_*Y\arrow CH_*(X\times Y)$$
is surjective for any open subset $Y\subset A^n$, since $X\times Y$
is an open subset of $X\times A^n$.

Now let $Y$ be any variety, $S\subset Y$ a closed subvariety,
and $U=Y-S$. Let $X$ be a variety. Then we have exact sequences as shown.
$$\begin{CD}
CH_*X\otimes_{\Z }CH_*S   @>>>   CH_*X\otimes_{\Z }CH_*Y @>>>   
   CH_*X\otimes_{\Z }CH_*U @>>> 0\\
@VVV                         @VVV    @VVV\\
CH_*(X\times S)  @>>>  CH_*(X\times Y) @>>> CH_*(X\times U) @>>> 0.
\end{CD}$$
It follows that if $CH_*X\otimes CH_*S\arrow CH_*(X\times S)$
and $CH_*X\otimes CH_*U\arrow CH_*(X\times U)$ are surjective,
then so is $CH_*X\otimes CH_*Y \arrow CH_*(X\times Y)$.
The theorem follows. 
\end{proof}

\begin{lemma}
\label{kunneth}
If $X$ and $Y$ are varieties such that the maps
$CH_*X\arrow H_*^{BM}(X,\Z )$ and $CH_*Y\arrow H_*^{BM}(Y,\Z )$
are split injections of abelian groups, then the map
$$CH_*X\otimes _{\Z }CH_*Y\arrow CH_*(X\times Y)$$
is injective. If in addition  $X$ or $Y$ can be partitioned into
open subsets of affine spaces, then we have
$$CH_*X\otimes _{\Z }CH_*Y\cong CH_*(X\times Y),$$
and $X\times Y$ also has a split injection $CH_*(X\times Y)
\inj H^{BM}_*(X\times Y,\Z ).$
\end{lemma}

\begin{proof}
Use the diagram
$$\begin{CD}
CH_*X\otimes _{\Z }CH_*Y @>>> CH_*(X\times Y)\\
@VVV                          @VVV\\
0\arrow H^{BM}_*X\otimes _{\Z }H^{BM}_*(Y,\Z )@>>> H^{BM}_*(X\times Y,\Z ),
\end{CD}$$
where the bottom arrow is split injective.
\end{proof}

Fulton has asked whether the Chow K\"{u}nneth formula might be true
assuming only that $Y$ can be cut into open subsets of affine spaces,
with no assumption on $X$ and no homological assumption.
In particular, Fulton, MacPherson, Sottile, and Sturmfels proved
this when $Y$ is a (possibly singular) toric variety or spherical variety
\cite{FMSS}. I generalized their result to an arbitrary ``linear variety''
$Y$ in the sense of Jannsen \cite{Totarolinear}.

In particular, let $G_1$ be any algebraic group such that there are quotients
$(V-S)/G_1$ with the codimension of $S$ arbitrarily large which
are linear varieties; this is true for most of the groups $G_1$ for which
we can compute the Chow ring of $BG_1$. Then
$$CH^*B(G_1\times G_2)=CH^*BG_1\otimes_{\Z}CH^*BG_2$$
for all algebraic groups $G_2$. Of course, this would not be true
for integral cohomology in place of the Chow ring, say for $G_1=\Z/p$.
But something similar is true for complex cobordism: if $G_1$ is a compact
Lie group such that the Morava $K$-theories of $BG_1$ are concentrated
in even degrees, then
$$MU^*B(G_1\times G_2)=MU^*BG_1\otimes_{MU^*}MU^*BG_2$$
for all compact Lie groups $G_2$ \cite{RWY}.

A special case of all this is that for all abelian groups $G$ over $\C$,
the Chow ring of $BG$ is the symmetric algebra on $CH^1BG=\text{Hom}(G,\C^*)$.
This ring maps isomorphically to $MU^*BG\otimes_{MU^*}\Z$,
by Landweber's calculation of $MU^*BG$ for $G$ abelian \cite{Landweberab}.

\section{Chow groups of cyclic products, introduction}
\label{CPintro}

We construct operations from the Chow groups of a quasi-projective
scheme $X$ to the Chow groups
of the $p$th cyclic product $Z^pX:=X^p/(\Z/p)$. We need $X$ quasi-projective
in order to know that this quotient variety $Z^pX$ exists; it will again
be quasi-projective.
From the construction
it will be clear that
these operations agree
with the analogous  operations on integral homology, as constructed
by Nakaoka \cite{Nakaoka56}. More precisely, Nakaoka constructed
such operations which raise degree by either an even or an odd
amount, and the operations we construct agree with his even operations.

\begin{lemma}
\label{ops}
Let $p$ be a prime number.
Let $X$ be a quasi-projective scheme over the complex numbers,
and let $Z^pX=X^p/(\Z /p)$ denote the $p$th cyclic product of $X$.
Then we will define operations $$\alpha ^j_i:CH_iX\arrow CH_jZ^pX, 
\mbox{  }i+1\leq j\leq pi-1$$ and
$$\gamma _i:CH_iX\arrow CH_{pi}Z^pX, \mbox{  }i>0,$$
with the following properties. (The formulas use the natural map
$(CH_*X)^{\otimes p}\arrow CH_*X^p\arrow CH_*Z^pX$.)

(1) $\alpha ^j_i$ is a homomorphism of abelian groups.

(2) $p\alpha^j_i=0$.

(3)$$\mbox{  }\gamma_i(x+y)=\gamma _ix+\sum\alpha_1\otimes\cdots\otimes\alpha_p
+\gamma_iy,$$
where the sum is over a set of representatives $\alpha =(\alpha_1,\cdots,
\alpha_p)$ for the action of $\Z /p$ on the set $\{x,y\}^p-\{(x,\cdots,
x),(y,\cdots,y)\}$.

(4) $p\gamma _ix=x^{\otimes p}$.
\end{lemma}

In the statement we assumed that $X$ was a scheme over the complex
numbers. The proof will in fact construct operations with all the same
properties over any field $k$ of characteristic $\neq p$ which
contains the $p$th roots of unity; the only complication is that
the operations $\alpha^j_i$ are naturally viewed as maps:
$$\alpha^j_i:CH_iX\otimes\mu_p(k)^{\otimes (pi-j)}\arrow CH_jZ^pX.$$

\begin{proof}
Let $C$ be any algebraic cycle of dimension $i>0$  on $X$.
(By definition, an algebraic cycle on $X$ is simply a $\Z$-linear combination
of irreducible subvarieties of $X$.) There is an associated cycle
$C^p$ on the product $X^p$. Since the subvarieties of $X$ occurring in $C$
have dimension $i>0$,  it is easy to see that
none of the subvarieties occurring in the cycle $C^p$ is
contained in the diagonal $\Delta_X\cong X\subset X^p$. Since the cycle
$C^p$ is invariant under the permutation action of the group $\Z/p$ on $X^p$,
the restriction of the cycle $C^p$ to $X^p-X$ (where $\Z/p$
acts freely) is the pullback of a unique cycle on $Z^pX-X$
under the etale morphism $X^p-X\arrow Z^pX-X$. Define
the cycle $Z^pC$ on $Z^pX$ to be the closure of this cycle on
$Z^pX-X$. Define the operation $\gamma_i$ on a cycle $C$
on $X$ to be the class of this cycle $Z^pC$ on $Z^pX$. We will show
below that this operation is well-defined on rational equivalence classes.

First we define the other operations $\alpha^j_i$. Let $C$ be an
algebraic cycle on $X$ of dimension $i>0$. Let $D$ denote the support
of $C$, that is, the union of the subvarieties of $X$ that occur in
the cycle $C$. There is an obvious
$\Z /p$-principal bundle over $Z^pD-D$.
Given a $p$th root of unity
and thus a homomorphism
$\Z /p \arrow k^*$,  we get a line bundle $L$ over $Z^pD-
\Delta _D$. (For $k=\C$, we use the $p$th root of unity $e^{2\pi i/p}$.)
Clearly $pc_1L=0\in \text{Pic}(Z^pD-\Delta _D )$. Line bundles act
on the Chow groups, lowering degree by 1, so we can
consider the element
$$c_1(L)^{pi-j}[Z^pC]$$
of $CH_{j}(Z^pD -\Delta _D)$ for $i+1\leq j\leq pi-1$.
Here $Z^pC$ is the cycle defined
in the previous paragraph; this cycle is supported
in the closed subset $Z^pD\subset Z^pX$, and the formula refers
to its restriction to a cycle on $Z^pD-D$.
The elements defined by this formula are 
killed by $p$ since $c_1(L)$ is. Since these classes are above the dimension
of $\Delta _D$, they uniquely determine elements
$$Y_j\in CH_jZ^pD,$$
$i+1\leq j\leq pi-1$. 
Clearly $pY_j=0$ in $CH_jZ^pD$ since this is true after restricting
to the open subset $Z^pD-D$. We define the operation
$\alpha^j_i$ on the cycle $C$
to be the image of this class $Y_j$ in $CH_iZ^pX$.

Let us show that the elements $\alpha^j_iC$ and $\gamma_iC$
in the Chow groups of
$Z^pX$ are well-defined when we replace the cycle $C$  by a 
rationally equivalent cycle. We have to show that for every
$(i+1)$-dimensional irreducible subvariety 
$W\subset X\times\Proj^1$ with the second projection not
constant, if we let $C_0$ and $C_{\infty}$ denote the
 fibers of $W$ over $0$ and $\infty$ in $\Proj^1$, viewed as  cycles on $X$,
then $\alpha^j_i(C_0)=\alpha^j_i(C_{\infty})$ and $\gamma_i(C_0)
=\gamma_i(C_{\infty})$ in the Chow groups of  $Z^pX$.
 
To see this, we use the fiber product $W^p/\Proj^1$. The group $\Z/p$
acts on this fiber product in a natural way, acting as the identity
on $\Proj^1$, and we call the  quotient scheme $Z^pW/\Proj^1$. The
fibers over 0 and infinity, as cycles on $Z^pX$, are exactly the cycles
$Z^pC_0=\gamma_i(C_0)$ and $Z ^pC_{\infty}=\gamma_i(C_{\infty})$.
Thus we have a rational equivalence
between the cycles $\gamma_i(C_0)$ and $\gamma_i(C_{\infty})$, so that
the operation $\gamma_i$ is well-defined on Chow groups.

We proceed to show that the other operations $\alpha^j_i$  are also
well-defined on Chow groups. Let $j$ be any integer with
$i+1\leq j\leq pi-1$.
As in the definition of $\alpha^j_i$, fix
 a $p$th root of unity in $k$ and thus a homomorphism $\Z/p\arrow k^*$.
This gives a   line bundle $L$ on $Z^pW/\Proj^1-W $.
Using the action of line bundles on Chow groups,
we get a $j+1$-dimensional Chow class on $Z^pW/\Proj^1-W$ associated
to the action of $c_1(L)^r$ on the fundamental cycle of the scheme
$Z^pW/\Proj^1-W$, where we let $r=pi-j$. Let $Y_j$ be the closure
in $Z^pW/\Proj^1$ of any cycle
on $Z^pW/\Proj^1-W$ representing this class.
The fibers of the cycle
$Y_j$ over $0$ and $\infty$ are $j$-dimensional
cycles on $Z^pW_0$ and $Z^pW_{\infty}$,
respectively, whose restrictions
 to $Z^pW_0-W_0$ and $Z^pW_{\infty}-W_{\infty}$
are rationally equivalent (on these open subsets)
to the cycles $\alpha^j_i(C_0)$ and $\alpha^j_i(C_{\infty})$. But
$W_0$ and $W_{\infty}$ have dimension $i$, which is less than
the dimension $j$ of these cycles. So in fact the fibers of $Y_j$
over $0$ and $\infty$ are rationally  equivalent to
the cycles $\alpha^j_i(C_0)$ and $\alpha^j_i(C_{\infty})$ on the whole sets
$Z^pW_0$ and $Z^pW_{\infty}$. By mapping $Y_j$ into $X$, it follows
that $\alpha^j_i(C_0)$ and $\alpha^j_i(C_{\infty})$ are rationally equivalent
on $X$, as we want.

We now prove properties (1)-(4) of the operations $\alpha ^j_i$ and
$\gamma_i$. It is easy to prove properties (1) and (3), describing
what happens to the operations after multiplying by $p$. Namely
$\gamma _i[C]$ is represented by the cycle $Z^pC$ on $Z^pX$, and $p$
times this cycle is the image of the cycle $C^p$ on $X^p$, since the map
$C^p\arrow Z^pC$ is generically $p$ to 1. (We are using that
$i=\mbox{dim }C$ is greater than 0 here.) And the classes
$\alpha^j_iC$ are obviously $p$-torsion, since the Chern class
$c_1L$ used to define them is $p$-torsion.

It is not much harder to describe what $\gamma _i$ of a sum of two cycles
is. We have
$$Z^p(C_1+C_2)=Z^pC_1+\sum f_*(\alpha _1\times\cdots\times\alpha_p)+Z^pC_2,$$
where $f$ denotes the map $X^p\arrow Z^pX$ and the sum runs over some
set of representatives $\alpha $ for the orbits of the group $\Z /p$
on the set $\{C_1,C_2\}^p-\{(C_1,\cdots,C_1), (C_2,\cdots,C_2)\}$.
This immediately implies the formula for $\gamma _i(x+y)$.

Finally, the formula above implies an equality in $CH_*(Z^pX-X )$:
$$c_1(L)^r[Z^p(C_1+C_2)]=c_1(L)^r[Z^pC_1]+\sum c_1(L)^rf_*[\alpha _1\times
\cdots\times\alpha_p]+c_1(L)^r[Z^pC_2],$$
for all $r\geq 0$. For $r\geq 1$, since $f^*c_1(L)=0$, the terms in the middle
are all 0, so we have that
$$\alpha^j_i(x+y)=\alpha^j_ix+\alpha^j_iy.$$
To be honest, we have only checked this equality
in $CH_*(Z^pX-X )$.
But we can apply this argument with $X$ replaced by its subscheme
$C_1\cup C_2$. Then we have proved that
$$\alpha^j_i(x+y)=\alpha^j_ix+\alpha^j_iy$$
in the Chow groups of $Z^p(C_1\cup C_2)-(C_1\cup C_2)$. Since these
elements of the Chow groups have dimension greater than that of
$C_1\cup C_2$, we have the same equality in the Chow groups
of $Z^p(C_1 \cup C_2)$. This implies the equality in the
Chow groups of $Z^pX$.
\end{proof}

We can sum up what we have done by defining a certain functor
$F_p$ from graded abelian groups $A_*$ to graded abelian groups;
Lemma \ref{ops} amounts to the assertion that for any variety $X$
there is a natural map $F_pCH_*X\arrow CH_*Z^pX$. Namely, let
$F_pA_*$ be the graded abelian group generated by $A_*\otimes_{\Z }\cdots
\otimes _{\Z }A_*$ ($p$ copies of $A_*$) together with
$\mu_p(k)^{\otimes (pi-j)}\otimes _{\Z }A_i$ in degree $j$ for
$i+1\leq j \leq pi-1$, and elements
$\gamma _ix_i$
in degree $pi$, where $x\in A_i$ and $i>0$.
To get $F_pA_*$ we divide out by the
relations
\begin{align*}
x_1\otimes\cdots\otimes x_p &= x_2\otimes \cdots\otimes x_p\otimes x_1\\
p\gamma_ix & =  x^{\otimes p}\\
\gamma_i(x+y)&= \gamma_ix+\sum\alpha_1\otimes \cdots\otimes \alpha_p+\gamma_i
y
\end{align*}
Here, as before, the sum in the last formula runs over the $\Z /p$-orbits
$\alpha $ in the set $\{x,y\}^p-{\{x,y\} }$.

\section{Chow groups of cyclic products, concluded}

We can now calculate the Chow groups of cyclic products for the
same class of varieties for which we proved the Chow K\"{u}nneth
formula (about ordinary products). Again we work over $\C$,
although the same proofs work for
varieties over any field $k$ of characteristic $\neq p$ which
contains the $p$th roots of unity, using etale homology
in place of ordinary homology. Actually, the following Lemma is stated
in terms of Borel-Moore homology, which can be defined as the homology
of the chain complex of locally finite singular chains on a locally
compact topological space \cite{Fulton}.

Since we cannot even
compute the Chow groups of a product in general, as discussed in
section \ref{Kunneth}, we cannot expect
to compute the Chow groups of arbitrary symmetric products. Fortunately
we only need to analyze symmetric products of a special class of
varieties in order to compute the
Chow ring of the symmetric group.

\begin{lemma}
\label{functor}

(1) There is a functor $F_p$, defined in section \ref{CPintro},
from graded abelian groups to graded abelian groups, such that there
is a natural map $F_pCH_*X\arrow CH_*Z^pX$ for any quasi-projective
variety $X$.

(2) If $X$ can be cut into open subsets of affine spaces, then so can
$Z^pX$, and the map $F_pCH_*X\arrow CH_*Z^pX$ is surjective.

(3) If $CH_*X\arrow H_*^{BM}(X,\Z )$ is split injective, then the
composition $F_pCH_*X\arrow CH_*Z^pX\arrow H_*^{BM}(X,\Z )$
is split injective.

Thus if $X$ satisfies hypotheses (2) and (3), then so does $Z^pX$,
and $F_pCH_*X\cong CH_*Z^pX$.
\end{lemma}

Here (1) was proved at the end of the last section.

\begin{proof} (2)
This
follows from the observations that the cyclic product $Z^pA^i$
can be cut into a union of open subsets of affine spaces of
dimensions $i+1,i+2,\dotsc, pi$, corresponding to the operations
$\alpha^j_i$ and $\gamma_i$ on the fundamental class of $A^n$,
and that an open subset of affine space has Chow groups equal to
$\Z$ in the top dimension and 0 below that.

To prove this description of $Z^pA^i$, we can analyze, more generally,
any quotient variety $V/(\Z/p)$, where $V$ is a representation of $\Z/p$
over a field $k$ in which $p$ is invertible and which contains
the $p$th roots of unity.
Suppose $V$ is a nontrivial representation of $\Z/p$. Then we can write
$V=W\oplus L$ where $L$ is a nontrivial 1-dimensional representation of $\Z/p$.
The quotient variety $V/(\Z/p)$ can be cut into $W/(\Z/p)$ and a $W$-bundle
over
$L^{\otimes p}-0\cong A^1-0$. The latter bundle can be seen, for example
by direct calculation, to be isomorphic to $W\times (A^1-0)$ as a variety.
This analysis gives by induction a decomposition of $V/(\Z/p)$ into
open subsets of affine spaces, as we need.
\end{proof}

\begin{proof} (3)
Suppose that
 $CH_*X\arrow H_*^{BM}(X,\Z )$ is split injective. Write the 
finitely generated abelian group $CH_*X$ as a direct sum of cyclic groups
$(\Z/a_i)\cdot e^i$, where $a_i$ is 0 or a prime power, $e^i\in CH_*X$,
and $1\leq i\leq n$. Let $S\subset \{1,\ldots,n\}$ be the set of $i$
such that $a_i$ is 0 or a power of $p$ and $\mbox{dim }e^i>0$.
Then, from the definition of the functor $F_p$, it is easy to check
that $F_pCH_*X$ is a quotient of the group
$$A_*:=\oplus \Z/(a_{i_1},\ldots,a_{i_p})\cdot e^{i_1}\otimes\cdots\otimes
e^{i_p}\oplus \oplus_i \Z/(pa_i)\cdot \gamma(e^i) \oplus \oplus_{i,j}
\Z/p\cdot\alpha^j(e^i).$$
Here the first sum runs over a set of orbit representatives for the action of
$\Z /p$ on $\{1,\ldots,n\}^p-S$, the second sum is over
$i\in S$, and the second sum is over $i\in S$ and $\mbox{dim }e^i+1\leq
j\leq p \mbox{ dim }e^i-1$.

We will show that the map $A_*\arrow H_*^{BM}(Z^pX,\Z )$ is split injective.
This will imply that $F_pCH_*X$ is actually isomorphic to $A_*$, not
just a quotient of it, and also that the composition
$F_pCH_*X\arrow CH_*Z^pX\arrow H_*^{BM}(Z^pX,\Z)$ is split injective,
thus proving (3).

We use the compatibility (which is clear) of our operations
$\gamma_i$ and $\alpha^j_i$ on Chow groups with the even-degree
operations from the homology of $X$ to the homology
of $Z^pX$ defined by Nakaoka \cite{Nakaoka56}. From Nakaoka's basis
for the homology of $Z^pX$ in terms of his operations,
we see that $F_p(H_*^{BM}(X,\Z))$ maps by a split injection
 into $H_*^{BM}(Z^pX,\Z)$.
Since $CH_*X$ is a direct summand of $H^{BM}_*(X,\Z)$, we deduce that
$A_*=F_pCH_*X$ maps by a split injection into  $H_*^{BM}(Z^pX,\Z)$.
\end{proof}

\section{Wreath products}
\label{detection}

\begin{lemma}
\label{wreath}
Let $G$ be a finite group which is an iterated wreath product
$\Z/p \wr \Z/p \wr \cdots \wr \Z/p$. Then the natural
homomorphism from the Chow ring of $(BG)_{\C}$ to the integer cohomology ring
of $BG$ is injective, and in fact additively split. The same is true
for products of groups of this form.
\end{lemma}

\begin{proof}
It
follows from Lemma \ref{functor} that if $G$ is a finite group
such that $BG$ can be approximated by smooth quasi-projective
varieties satisfying
conditions (2) and (3) in the lemma, then the wreath product
$\Z /p \wr G$ has the same property. That is, writing out conditions
(2) and (3) in more detail, the wreath product $\Z /p \wr G$
has the property that $B(\Z /p \wr G)$ can be approximated by
smooth varieties which can be cut into open subsets of affine spaces,
and the map $CH^*B(\Z/p \wr G) \arrow H^*(B(\Z /p \wr G),\Z)$
is split injective. In deducing this statement from Lemma
\ref{functor}, the point is that if $X$ is a
smooth variety which approximates
$BG$ up to some high dimension, then the complement of the diagonals
in the cyclic product $Z^pX$ is a smooth variety which approximates
$B(\Z/p \wr G)$, and those diagonals have high codimension,
so that the Chow cohomology in low degrees of $Z^pX-(\mbox
{diagonals})$ is equal to the low-codimension Chow groups of 
the singular variety $Z^pX$;
and those groups are what Lemma \ref{functor} computes.

This implies the lemma for wreath products of copies of $\Z/p$,
by induction. For products of such groups, the lemma follows,
using the Chow K\"{u}nneth formula (Lemma \ref{kunneth}).
\end{proof}

Lemma \ref{functor} actually gives an explicit calculation of the
Chow cohomology groups of $BG$,
for groups $G$ as in Lemma \ref{wreath}. In particular, we can see from
this calculation that for such groups $G$, the mod $p$ Chow ring
of the wreath product $\Z /p \wr G$ is detected on the two subgroups
$G^p$ and $\Z /p \times G$. This means that the homomorphism
$$CH^*B(\Z /p \wr G)/p \arrow CH^*(BG^p)/p \oplus CH^*B(\Z/p\times G)/p$$
is injective. Applying this repeatedly, we find that
 the mod $p$ Chow ring of a group $G$ as in Lemma \ref{wreath}
is detected on the elementary abelian subgroups of $G$.

\section{The symmetric group: injectivity}
\label{inj}

The calculation of $CH^*BS _n$, locally at the prime $p$ (that is,
the calculation of $CH^*BS _n \otimes_{\Z }\Z_{(p)}$),
follows from the calculation of $CH^*B(\Z /p\wr\cdots\wr\Z /p)$.
Namely, let $H$ be the $p$-Sylow subgroup of the symmetric group
$S _n$; if one writes $n=p^{i_1}+p^{i_2}+\cdots$, $i_1<i_2<\cdots $,
then $H=(\Z /p)^{\wr i_1}\times (\Z /p)^{\wr i_2}\times\cdots $,
where $(\Z/p)^{\wr i}$ denotes the $i$-fold wreath product
$\Z/p \wr\cdots \wr \Z/p$.
 We know
the Chow cohomology groups of $BH$ by Lemma \ref{wreath}
and the comments afterward.
In particular we know that the Chow cohomology
groups of $BH$ inject onto a direct summand of its integer cohomology.
The $p$-localization of the Chow groups of $BS _n$, i.e.,
$CH^*(BS _n)_{(p)}:=CH^*(BS  _n)\otimes _{\Z }\Z _{(p)}$,
can then be described  as follows.

\begin{lemma}
$$CH^*(BS _n)_{(p)} =CH^*BH\cap H^*(BS _n,\Z )_{(p)}\subset
H^*(BH,\Z )_{(p)}.$$
\end{lemma}

In particular, the Chow ring of the symmetric group maps
injectively to its integral cohomology, by applying this lemma for
each prime number $p$.
In fact the proof 
will show that this map is a split injection.

\begin{proof}
A standard transfer argument shows that if $H$ denotes
the $p$-Sylow subgroup of a finite group $G$, then we have a split
injection
$$H^*(BG,\Z)_{(p)}\inj H^*(BH,\Z),$$
and the image is precisely the ``$G$-invariant'' subgroup of $H^*BH$
(\cite{Brown}, p.~84).
Here we say that  $x\in H^*BH$ is $G$-invariant if $x|_{H\cap
gHg^{-1}}=gxg^{-1}|_{H\cap gHg^{-1}}$ for all $g\in G$. The same
argument gives an injection of Chow rings:
$$CH^*(BG)_{(p)}\inj CH^*BH$$
is the subgroup of $G$-invariant elements of $CH^*BH$.

Thus the inclusion $CH^*(BG)_{(p)}\subset CH^*BH\cap H^*BG_{(p)}$
is obvious. To prove the opposite direction, we have to show that
an element of $CH^*BH$ whose image in $H^*BH$ is $G$-invariant is
itself $G$-invariant. This is not obvious, since we do not know whether
$CH^*B(H\cap gHg^{-1})$ injects into $H^*B(H\cap gHg^{-1})$
even though we know this for $H$ itself. But in fact it is easy.
Let $x\in CH^*BH$ have $G$-invariant image in $H^*BH$. Then the standard
transfer argument, already used above, shows that
$\text{res}^G_H \text{tr}^G_Hx=(G:H) x\in H^*BH$. But the restriction maps
in $CH^*$ and $H^*$ are compatible, and likewise for transfer.
And $CH^*BH$ injects into $H^*BH$, so we actually have $\text{res}^G_H
 \text{tr}^G_H
x=(G:H)x\in CH^*BH$. Since $(G:H)$ is a $p$-local unit, this implies that $x$
is the restriction of an element of $CH^*BG$ and so is $G$-invariant.
\end{proof}

The same transfer argument shows that $CH^*(BS _n)\otimes\Z/p$
injects into $H^*(BS _n,\Z/p)$. It follows that
the map $CH^*(BS _n)_{(p)}\arrow H^*(BS _n,\Z)_{(p)}$
is a split injection.

\section{The Chow ring of the symmetric group}
\label{trans}

\begin{theorem}
\label{symm}
The homomorphism
$$CH^*BS_n \arrow MU^*BS_n\otimes_{MU^*}\Z$$
 is an isomorphism. These rings map by an additively split injection
to the cohomology ring $H^*(BS_n,\Z)$.
\end{theorem}

\begin{proof}
In section 9, we showed that the Chow ring of
$BS_n$ maps additively split injectively to  its cohomology.
A fortiori,
the Chow ring injects into $MU^*BS_n\otimes_{MU^*}\Z$.
This injection
 is in fact an isomorphism, since Hopkins-Kuhn-Ravenel \cite{HKR}
and Ravenel-Wilson-Yagita \cite{RWY} have shown that the topological
$MU^*$-module  $MU^*BS_n$
is generated by transferred Euler classes, and the image of such a class
in the quotient ring $MU^*BS_n\otimes_{MU^*}\Z$ comes from $CH^*BS_n$.
\end{proof}

This proof shows that the Chow ring of any symmetric group
is generated by transferred Euler classes. That can also be seen more
directly by making our computation of
the Chow ring of wreath products more explicit. Indeed, the Chow ring
of iterated wreath products of $\Z/p$ is generated by transferred
Euler classes, and this implies the same statement for the
symmetric group.

\section{The ring $CH^*(BS _n)\otimes \Z/2$}
\label{mod2}

We now describe the ring $CH^*(BS _n)\otimes \Z/2$ more explicitly.

\begin{proposition}
\label{calc2}
For the symmetric group $S _n$, the two maps
$$\mbox{Frob}: H^i(BS_n,\Z /2)\arrow H^{2i}(BS_n,\Z /2)$$
($x\mapsto x^2$) and 
$$CH^i(BS_n)\otimes \Z/2\arrow H^{2i}(BS_n,\Z /2)$$
are injective ring homomorphisms with the same image. This implies that
we have a ring isomorphism between $H^*(BS_n,\Z /2)$ and 
$CH^*(BS_n)\otimes \Z /2$,
with  $H^i$ corresponding to $CH^i$.
\end{proposition}

\begin{proof}
The first homomorphism is injective
because the ring $H^*(BS_n,\Z /2)$ has no nilpotents; this in turn follows
from the $\Z /2$-cohomology of these groups being detected by elementary
abelian subgroups (\cite{Madsen-Milgram}, pp.~53-57). The second
homomorphism is injective by Theorem \ref{symm}.

Next, by the calculation of the $\Z/2$-cohomology of wreath products
of copies of $\Z/2$, $H^*(BS _n,\Z /2)$ is generated by transferred
Euler-Stiefel-Whitney classes, that is, transfers of
top Stiefel-Whitney classes
of real representations of subgroups of $S _n$. Then observe that
the squaring map in $\Z /2$-cohomology,
$$\mbox{Frob}:H^i(X,\Z/2)\arrow H^{2i}(X,\Z /2),$$
coincides for each $i$ with a certain Steenrod operation, thus a stable
operation, which therefore commutes with transfers for finite coverings.
That is, for a finite covering $f:X\arrow Y$ and $x\in H^*(X,\Z /2)$,
$f_*(x^2)=(f_*x)^2$. It follows that $H^*(BS _n,\Z /2)^2$ is generated by
transferred Euler classes of complex representations of subgroups,
since for a real representation $E$ of a group
$H$, $\chi(E\otimes \C)=\chi(E)^2\in
H^*(BH,\Z /2)$. Thus we have the inclusion
$$H^*(BS _n,\Z /2)^2\subset \im(CH^*(BS _n)\otimes \Z/2\arrow
H^*(BS _n,\Z /2)).$$

To get the inclusion the other way, it suffices to check that
$CH^*(BS _n)\otimes \Z /2$ is generated by transferred Euler classes
of complexified real representations. We showed in section \ref{trans}
that for any prime $p$, $CH^*(BS _n)\otimes \Z /p$ is generated
by transferred Euler classes of complex representations. These representations
are obtained by an inductive procedure starting from the basic 
representation $\Z /p\arrow \C ^*$. But for $p=2$, this representation
is the complexification of a real representation, and so the same
applies to all the resulting representations.
\end{proof}

This strange isomorphism has the property that the Stiefel-Whitney
classes of the standard representation $S_n\arrow GL(n,\R)$
correspond to the Chern classes of the standard representation
$S_n\arrow GL(n,\C)$. Just as the ring $CH^*BS _6$ is not generated
by Chern  classes, the ring $H^*(BS _6,\Z/2)$ is not generated
by Stiefel-Whitney classes.

\section{The Chow ring of the symmetric group over a general field}
\label{fields}

Let $G$ be a finite group.
We can view $BG$ as a limit of varieties defined over $\Q $, or even
over $\Z $. By base change we get $(BG)_k$ for any field $k$.
So far we have computed $CH^*(BG)_{k,(p)}$ when $G$ is the symmetric group
and $k$ is a field of characteristic $\neq p$ which contains the $p$th
roots of unity. In this section we show that the Chow ring $\otimes \Z_{(p)}$
of the
symmetric group is the same for all
fields $k$ of characteristic $\neq p$.

This result can be divided into two parts. First, for any finite group $G$
there is an action of $Gal(\overline{k}/k)$ on $CH^*(BG)_{\k}$,
and for some groups such as the group of order $p$ and the symmetric group,
we find that the natural map
$$CH^*(BG)_k\arrow CH^*(BG)_{\k}^{Gal(\k /k)}$$
is an isomorphism. Then for the symmetric group, but not for the group
of order $p$ when $p\geq 3$, we show that $Gal(\k /k)$ acts trivially on $CH^*BG$.

Even the first statement
is not true for arbitrary finite groups. Even for $K$-theory,
which is simpler than the Chow groups, and even if the
characteristic of $k$ is prime to the order of $G$, the map
$$K_0BG_k\arrow (K_0BG_{\k})^{Gal(\k /k)}$$
is not always an isomorphism, although it is an isomorphism
modulo torsion. The map fails to be an isomorphism if and only if
$G$ has a representation over $k$ with Schur index not equal to 1
\cite{Serre}, for example when $G$ is the quaternion group and $k=\R $.

When $G$ is the symmetric group, all Schur indices are equal to 1,
and in fact all representations of $G$ are defined over $\Q $. (We will
not use these facts in what follows.) In particular, the action of
$Gal(\overline{\Q}/\Q)$ on $K_0BS _n$ is trivial. It follows that
the action of $Gal(\overline{\Q}/\Q)$ on $CH^*BS _n$ fixes all Chern classes
of representations of $S _n$; but we need a different argument,
below, to show that the Galois action on all of
$CH^*BS _n$ is trivial.

For any finite group $G$, we can analyze the Galois action on $CH^*(BG)_{(p)}$
using the cycle map
$$CH^i(BG)_{\k ,(p)}\arrow H^{2i}(BG_{\k},\Z_p(i)),$$
which is Galois-equivariant. For the groups $G$ we consider, this homomorphism
is injective, so it suffices to describe the Galois action on the group
$H^{2i}(BG_{\k},\Z_p(i))=H^{2i}(BG_{\k},\Z_p)\otimes_{\Z_p}\Z_p(i)$.
The answer is simple: for all finite groups $G$, $Gal(\k /k)$ acts
trivially on $H^*(BG,\Z _p)$. This is more or less obvious from the 
definition of etale cohomology, since we can view $BG$ as a limit
of quotient varieties $X/G$ in which $X$ is defined over $k$ and
every element of $G$ acts by a map of varieties over $k$.

Thus the Galois action on $H^{2i}(BG_{\k },\Z_p(i))$ just comes from
the Galois action on $\Z_p(i)$. In the extreme case where $k=\Q $,
we can describe the Galois-fixed subgroup of $H^{2i}(BG,\Z_p(i))$
as follows, by a simple calculation given by Grothendieck \cite{Grothendieck}:
\begin{equation*}
 H^{2i}(BG,\Z_p(i))^{Gal(\overline{\Q}/\Q)}=
\begin{cases}
0 & \text{if $i\not\equiv 0\pmod{p-1}$}\\
\text{ker }p & \text{if }i=a(p-1), (a,p)=1\\
\text{ker }p^{r+1} & \text{if }i=ap^r(p-1), (a,p)=1
\end{cases}
\end{equation*}
The subgroup fixed by $Gal(\QQ /\Q(\mu_p))$ is given by the same restrictions
except that $i$ is not required to be a multiple of $p-1$.

\begin{example}
The cyclic group of order $p$ gives an interesting example
of the Galois action on $CH^*BG$. Namely, let $\overline{\Q}$
be the algebraic closure of $\Q$. An element of $CH^1(B\Z/p)_{\overline{\Q}}$
 is given
by a homomorphism $\Z/p\arrow \overline{\Q}^*$,
 so that $CH^1(B\Z/p)_{\overline{\Q}}
\cong \mu_p(\overline{\Q})$, and this description is natural with respect
to automorphisms of $\overline{\Q}$. Since the  ring $CH^*(B\Z/p)_{\QQ}$
is generated by $CH^1$, we have a natural isomorphism 
$CH^i(B\Z/p)_{\overline{\Q}}\cong (\mu_p(\overline{\Q}))^{\otimes i}$.
It follows that if we view $B\Z/p$ as a limit of varieties over $\Q$, say,
then the Chow ring of $(B\Z/p)_{\Q}$ maps to 0 except in dimensions
$\equiv 0\pmod{p-1}$. The point is that there is a natural isomorphism
$(\mu_p)^{\otimes (p-1)}\cong \Z/p$, because there is a natural generator
for the cyclic group $(\mu_p)^{\otimes (p-1)}$: take the $(p-1)$st power
of any generator of the cyclic group $\mu _p$. This argument gives
a canonical generator for $(\mu _p)^{\otimes i}$ 
whenever $i\equiv 0\pmod{p-1}$,
but otherwise there is none.

As a matter of fact the Chow ring of $(B\Z/p)_{\Q }$  is exactly the subring
of elements of dimension $\equiv 0\pmod{p-1}$, that is, a $\Z/p$-polynomial
algebra on one generator of degree $p-1$. As it happens this is also
the image of the Chow ring of $BS _p$ in the Chow ring of its
$p$-Sylow subgroup, $\Z/p$.
\end{example}

We know that the Chow ring $\otimes \Z_{(p)}$
of the symmetric group over $\Q(\mu _p)$
is the same as over $\QQ$. But the extension $\Q(\mu_p)/\Q$
has degree $p-1$ which is prime to $p$. So the usual transfer
argument for the etale covering
$$\begin{CD}
(BS _n)_{\Q(\mu_p)}\\
@V(\Z /p)^*VV \\
(BS _n)_{\Q }
\end{CD}$$
shows that $CH^*(BS _n)_{\Q, (p)}\cong CH^*(BS _n)_{\Q(\mu _p),
(p)}^{(\Z /p)^*}$. So we just have to show that the action of
$Gal(\QQ / \Q)$ through its quotient $(\Z /p)^*$ on $CH^*(BS _n)
_{\Q(\mu_p),(p)}$ is trivial.

But our earlier comments on the Galois action on $H^{2i}(BG,
\Z_p(i))$,  together with the injectivity of $CH^iBS _n\arrow
H^{2i}(BG,\Z_p(i))$, shows that $(\Z /p)^*$ acts trivially on
$CH^*(BS _n)_{(p)}$ if these $p$-local
Chow groups are concentrated in dimensions
$\equiv 0\pmod{p-1}$.

To prove this, we use that $CH^*(BS _n)/p$ is detected
on certain explicitly known elementary abelian subgroups of
$S _n$. Indeed, this mod $p$ Chow ring injects into the mod $p$ Chow ring
of the $p$-Sylow subgroup of $S_n$, which is a product
of iterated wreath products of copies of $\Z/p$, and we proved
in section \ref{detection} that, for a class of groups $G$ including
those considered here, the mod $p$ Chow ring of the wreath product of
$\Z /p \wr G$ is detected on the two subgroups $G^p$ and $\Z /p \times
G$.

Then we observe that for all of these elementary abelian
subgroups $H=(\Z/p)^k$ of $S _n$, the normalizer in $S_n$ contains
a group $(\Z/p)^*$ which acts in the obvious way (by scalar
multiplication) on $H$. It follows that $(\Z/p)^*$ acts
in the natural way on the polynomial ring $CH^*BH=
\Z/p[x_1,\ldots ,x_k]$, so that the invariants of $(\Z/p)^*$
in this polynomial ring is precisely the subring consisting
of elements of dimension a multiple of $p-1$. But the homomorphism
of the Chow ring of $S _n$ into that of $H$ automatically
maps into this ring of invariants. Combining this with the previous
paragraph, we have proved that the mod $p$ Chow ring of $S _n$
is nonzero only in dimensions a multiple of $p-1$. It follows that
the $p$-local Chow ring of $S_n$ is also concentrated in these
dimensions.

This is what we needed, three paragraphs ago, to complete the proof
that the Chow groups of the symmetric
group are the  same over all fields of characteristic 0.

We only need a few more words to explain why the $p$-local Chow groups
of the symmetric group are the same over any field of characteristic
$\neq p$.  Suppose $l$ is a prime number not equal to $p$. We have
shown that the Chow groups of $BS _n$ are the same over any field
which contains $F_l(\mu_p)$, and in fact the same as in characteristic 0.
 The Galois group $Gal(F_l(\mu_p)/F_l)$
is a subgroup of $(\Z/p)^*$, and since all of $(\Z /p)^*$ acts trivially
on the Chow groups of $S _n$, so does this subgroup. 
(All we are using here is that the $p$-local Chow groups of
$(BS_n)_{F_l(\mu_p)}$ are concentrated in dimensions
$\equiv 0\pmod{p-1}$.)
We deduce that
the Chow groups $CH^*(BS _n)_{(p)}$ are the same over all fields of
characteristic $\neq p$.

\section{Generators for the Chow ring of $BG$}
\label{bound}

In this section, we give a simple upper bound for the degrees of
a set of generators for the Chow ring of $BG$, for any
algebraic group $G$ (not necessarily connected). As mentioned
in the introduction, nothing similar is known for the ordinary
cohomology or the complex cobordism of $BG$.

\begin{theorem}
\label{gensthm}
Let $G$ be an algebraic group over a field, and let $G\inj H$
be an imbedding of $G$ into a group $H$ which is a product of the groups
$GL(n)$ for some integers $n$.
Then the Chow ring $CH^*BG$ is generated
as a module over $CH^*BH$ by elements of degree
at most $\text{dim }H/G$. Here $CH^*BH$ is just a polynomial ring
over $\Z$ generated by the Chern classes (see section \ref{classical}).
\end{theorem}

This follows from the more precise statement:

\begin{proposition}
\label{gens}
In the situation of the theorem, we have
$$CH^*(H/G)=CH^*BG\otimes_{CH^*BH}\Z.$$
\end{proposition}

Indeed, Proposition \ref{gens} implies that the Chow ring of $BG$ maps onto
that of $H/G$, and that $CH^*BG$ is generated as a $CH^*BH$-module
by any set of elements of $CH^*BG$ which restrict to generators
for the Chow ring $CH^*(H/G)$. So, in particular, $CH^*BG$ is generated
as a $CH^*BH$-module by elements of degree at most
$\text{dim }H/G$, thus proving Theorem \ref{gensthm}.  Proposition \ref{gens}
says, more generally, that to find generators for the Chow ring
of $BG$, it suffices to find generators for the Chow ring
of a single quotient variety $GL(n)/G$; this will be applied
in the next section to compute the Chow ring of certain
classifying spaces.

\begin{proof} (Proposition \ref{gens})
The point is that there is a
fibration
$$\xymatrix{H/G \ar[r]& BG\ar[d] \\
     & BH}$$
with structure group $H$. (Start with the universal fibration
$H\arrow EH\arrow BH$, and then form the quotient
$(EH\times H/G)/H$; this fibers over $BH$ with fiber $H/G$,
and it can be identified with $BG$.) To avoid talking about anything
infinite-dimensional, it suffices to consider the corresponding
fibration
$$\xymatrix{H/G \ar[r]& (V-S)/G\ar[d] \\
     & (V-S)/H}$$
associated to a representation $V$ of $H$.
Later in the proof we will restrict ourselves to a special class
of representations of $H$, but since these include representations
with $S$ of arbitrarily large codimension, that will be enough
to prove the proposition.

Since $H$ is on Grothendieck's list of ``special''
groups, every principal $H$-bundle in algebraic geometry
is Zariski-locally trivial \cite{Sem}. So the above fibration,
being associated to the principal $H$-bundle over $(V-S)/H$,
is  Zariski-locally trivial.

This implies that the restriction map $CH^*(V-S)/G\arrow
CH^*(H/G)$ associated to the fiber over a point $x\in(V-S)/H$ is surjective.
Indeed, choose a trivialization of the above bundle over a Zariski
open neighborhood $U$ of $x$ in $(V-S)/H$. Then, for any subvariety
$Z$ of the fiber over $x$, we can spread it out to a variety
$Z\times U$, and take the closure to get a subvariety of $(V-S)/G$.
This gives a subvariety  of $(V-S)/G$ which restricts to the given
subvariety $Z$ of the fiber $H/G$, so that $CH^*(V-S)/G\arrow
CH^*(H/G)$ is surjective.
Clearly elements of $CH^*(V-S)/H$
of positive degree restrict to 0 in the fiber $H/G$,
so we have a surjection 
$$CH^*(V-S)/G\otimes_{CH^*(V-S)/H}\Z\surj CH^*(H/G).$$

The proposition follows if we can prove that this map is
an isomorphism.
To see this, suppose that $x\in CH^*(V-S)/G$
restricts to 0 in the chosen fiber $H/G$; we have to show that
$x$ can be written as a finite sum $x=\sum a_ix_i$ with
$a_i\in CH^{>0}(V-S)/H$, $x_i\in CH^*(V-S)/G$. For this it seems natural
to use the detailed information which we possess about $BH$
since $H$ is a product of the groups $GL(n)$: namely, taking
$V$ to be the direct sum of copies of the standard representations of
the various factors $GL(n)$, we can arrange for $(V-S)/H$ to
be a product of Grassmannians, while making the codimension of $S$
as large as we like. In particular, using such representations $V$ of $H$,
$(V-S)/H$ is a smooth projective variety
which has a  cell decomposition, and the natural principal $H$-bundle
over $(V-S)/H$ is trivial over each cell. So also our $H/G$-bundle
over $(V-S)/H$ is trivial over each cell, and for each cell $C\subset BH$,
$(V-S)/G|_C$ is isomorphic to $(H/G)\times C$ and the Chow ring of
$(V-S)/G|_C$ is isomorphic to that of $H/G$. In particular,
an element $x\in CH^*(V-S)/G$ which restricts to 0 in
the Chow ring of a fixed fiber also restricts to
0 in the Chow ring of $(V-S)/G|_{C}$, where $C$ is the big cell
of $(V-S)/H$.

We now use the exact sequence
$$CH_*Y\arrow CH_*X\arrow CH_*(X-Y)\arrow 0$$
which exists for any algebraic variety $X$ and closed subset
$Y$. It follows that the above element $x$ is equivalent in 
$CH^*(V-S)/G$ to a linear combination of Zariski closures of product
varieties $Z\times C$ for cells $C$ in $(V-S)/H$ of codimension greater
than 0. By induction on the dimension of these cells, and using
that $CH^*(V-S)/G$ maps onto the Chow ring of every fiber $H/G$,
we find that $x$ can be written as a sum $\sum a_ix_i$ with
$a_i\in CH^{>0}(V-S)/H$ (representing the various cells of $(V-S)/H$
of codimension $>0$) and $x_i\in CH^*(V-S)/G$.
\end{proof}

\section{The Chow rings of the classical groups}
\label{classical}

For $G=GL(n,\C)$, then as we have
explained in the proof of Proposition \ref{gens},
 $BG$ can be approximated by smooth projective
varieties $(V-S)/G$
which have (algebraic) cell decompositions. For such varieties,
it is well known that the Chow ring maps isomorphically to the
cohomology ring \cite{Fulton}, and hence we have $CH^*BG=
H^*(BG,\Z)$ for such groups. In particular:

$$CH^*B\C^*=\Z[c_1]$$
$$CH^*BGL(n,\C)=\Z[c_1,\ldots,c_n]$$
Here and in what follows, we write $c_i$ for the $i$th Chern class
of the standard representation of a group $G$ when this has
an obvious meaning. These rings agree with $MU^*BG\otimes_{MU^*}\Z$,
for the trivial reason that $MU^*X\otimes_{MU^*}\Z=H^*(X,\Z)$
for all spaces $X$ with torsion-free cohomology.

A much more interesting calculation, proved below, is that
$$CH^*BO(n)=\Z[c_1,\ldots,c_n]/(2c_i=0\text{ for }i\text{ odd}).$$
This agrees with $MU^*BO(n)\otimes_{MU^*}\Z$ by Wilson's 
calculation of $MU^*BO(n)$ \cite{WilsonO}. It is simpler
than the integral cohomology of $BO(n)$, which has more
2-torsion related to Stiefel-Whitney classes. Let us mention
that $CH^*BO(n)\arrow H^*(BO(n),\Z)$ is injective,
and in fact additively split injective. This is
clear after tensoring  with $\Z_{(p)}$ for an odd prime $p$,
where both rings are polynomial rings on the Pontrjagin classes
$c_2,c_4,c_6,\ldots,c_{2\lfloor n/2 \rfloor}$.
To prove that the homomorphism
$$\Z[c_1,\ldots,c_n]/(2c_i=0\text{ for }i\text{ odd})\arrow H^*(BO(n),\Z)$$
 is additively split injective
$2$-locally, it suffices to show that it is injective after tensoring
with $\Z/2$, given that the first groups have no 4-torsion. 
To prove the injectivity after tensoring with $\Z/2$, we have
$$CH^*BO(n)\otimes_{\Z}\Z/2 = \Z/2[c_1,c_2,\ldots,c_n]$$
$$H^*(BO(n),\Z)\otimes_{\Z}\Z/2\inj H^*(BO(n),\Z/2)=\Z/2[w_1,w_2,
\ldots, w_n],$$
and $c_i$ restricts to $w_i^2$ \cite{MS}, so the map is injective.

To make the above computation of the Chow ring of $BO(n)$,
we use Proposition \ref{gens}, applied to the
standard representation $O(n)\subset GL(n)$. We find that
the Chow ring of $BO(n)$ is generated as a module over
$CH^*BGL(n)=\Z[c_1,\ldots, c_n]$ by any elements of $CH^*BO(n)$
which map onto $CH^*GL(n)/O(n)$. But the variety $GL(n)/O(n)$
is precisely the space of nondegenerate quadratic forms on
$\C^n$, and thus it is a Zariski open subset of affine space
$\C^{n(n+1)/2}$. By the fundamental exact sequence for Chow
groups, the Chow groups of any Zariski
open subset of affine space are 0 in codimension $>0$ (and $\Z$
in codimension 0). It follows that $CH^*BGL(n)=\Z[c_1,\ldots,
c_n]$ maps onto $CH^*BO(n)$.

It is easy to check that the relations $2c_i=0$ for $i$ odd are true
in $CH^*BO(n)$: this is because the representation $O(n)\arrow
GL(n)$ is self-dual, and we have $c_i(E^*)=(-1)^ic_i(E)$ in $CH^*X$
for every algebraic vector bundle $E$ on a variety $X$.
(This is standard
for algebraic geometers; for topologists, one might say that
it follows from the corresponding equation in $H^*(X,\Z)$ by considering
the universal bundle on $BGL(n)$ and using that $CH^*BGL(n)=
H^*(BGL(n),\Z)$.)

To show that the resulting surjection
$$\Z[c_1,\ldots,c_n]/(2c_i=0 \text{ for }i\text{ odd})\arrow CH^*BO(n)$$
is an isomorphism, it suffices to observe that, by the calculation above,
the composition $\Z[c_1,\ldots,c_n]/(2c_i=0\text{ for }i
\text{ odd})\arrow CH^*BO(n)\arrow H^*(BO(n),\Z)$ is injective,
so the first map is also injective. This proves the description
of $CH^*BO(n)$ stated above.

The symplectic group $Sp(2n)$ is similar but simpler. As in the case
of the orthogonal group, the quotient variety $GL(2n)/Sp(2n)$
is the space of nondegenerate alternating forms on $A^{2n}$
and hence is an open subspace of affine space. Thus $GL(2n)/Sp(2n)$
has trivial Chow groups, and by Proposition \ref{gens}, the homomorphism
$$CH^*BGL(2n)\arrow CH^*BSp(2n)$$
is surjective. Thus the Chow ring of $BSp(2n)$ is generated
by the Chern classes of the natural representation, $c_1,\ldots,c_{2n}$.

Since the natural representation of $Sp(2n)$ is self-dual, we have
$2c_i=0$ for $i$ odd in $CH^i(BSp(2n))$. In fact, we have $c_i=0$
for $i$ odd. To see this, it is enough to show that the Chow ring
of $BSp(2n)$ injects into the Chow ring of $BT$ for a maximal torus $T$,
since the latter ring is torsion-free. Here the classifying space
of a maximal torus can be viewed as an iterated affine-space bundle over the
classifying space of a Borel subgroup $B$, so it is equivalent
to show that the Chow ring of $BSp(2n)$ injects into the Chow ring
of $BB$. But $BB$ is a bundle over $BSp(2n)$ with fibers the smooth
projective variety $Sp(2n)/B$. Moreover, the group $Sp(2n)$ (unlike
the orthogonal group) is special \cite{Sem}, so that this bundle is
Zariski-locally trivial. Hence, taking the closure of a section of
this bundle over a Zariski open subset, there is an element
$\alpha\in CH^rBB$, where $r=\text{dim }Sp(2n)/B$, such that
$f_*\alpha=1\in CH^0BSp(2n)$. This implies that
$f^*:CH^*BSp(2n)\arrow CH^*BB$ is injective as we want, since
$$f_*(\alpha \cdot f^*x)=x$$
for all $x\in CH^*BSp(2n)$.

Thus the Chow ring of $BSp(2n)$ is generated by the even Chern classes
$c_2,c_4,\dotsc,c_{2n}$ of the standard representation. Using, for example,
the map from the Chow ring to cohomology, which is the polynomial ring
on these classes, we deduce that
$$CH^*BSp(2n)=H^*(BSp(2n),\Z)=\Z[c_2,c_4,\dotsc,c_{2n}].$$

\section{Other connected groups}
\label{other}

It is easy to compute the Chow ring of $BSO(2n+1)$ using the
previous  section's computation for $BO(2n+1)$. Indeed, there
is an isomorphism of groups
$$O(2n+1)\cong \Z/2\times SO(2n+1),$$
where the $\Z/2$ is generated by $-1\in O(2n+1)$. Since $B\Z/2$
can be approximated by linear varieties in the sense of
\cite{Totarolinear}, we have
$$CH^*B(G\times \Z/2)=CH^*BG\otimes_{\Z}CH^*B\Z/2,$$
where $CH^*B\Z/2=\Z[x]/(2x=0)$, $x\in CH^1$. Applying this
to $G=SO(2n+1)$ and using the previous section's computation
of $CH^*BO(2n+1)$, we find that
$$CH^*BSO(2n+1)=\Z[c_2,c_3,\ldots,c_{2n+1}]/(2c_i=0 \text{ for }
i\text{ odd}).$$
Since we also have $MU^*B(G\times \Z/2)=MU^*BG\otimes_{MU^*} B\Z/2$
for all compact Lie groups $G$, Wilson's calculation 
of $MU^*BO(2n+1)$ determines the calculation of $MU^*BSO(2n+1)$,
and we find that $CH^*BSO(2n+1)$ agrees with
$MU^*BSO(2n+1)\otimes_{MU^*}\Z$. Pandharipande has another computation
of the Chow rings of $BO(n)$ and $BSO(2n+1)$, and he used these
calculations to compute the Chow ring of the variety of rational
normal curves in $\Proj^n$ \cite{Pandharipande}.

The groups $SO(2n)$ are more complicated, in that the Chow ring
is not generated by the Chern classes of the standard representation.
One also has, at least, the $n$th Chern class of the representation
of $SO(2n)$
whose highest weight is twice that of one of the two spin representations
of $Spin(2n)$: this class is linearly independent of the space of polynomials
in the Chern classes
of the standard representation even in rational cohomology.
One can guess that the Chow ring of $SO(2n)$ is generated by
the Chern classes of the standard representation together
with this one further class; Pandharipande has proved
this in the case of $SO(4)$ \cite{Pandharipande}.
His calculation agrees with $MU^*BSO(4)\otimes_{MU^*}\Z$,
where $MU^*BSO(4)$ was computed by Kono and Yagita \cite{KY}.

The group $PGL(3)$ is even more complicated, in that
$MU^*BPGL(3)$ is concentrated in even degrees but
not generated by Chern classes,
by Kono and Yagita \cite{KY}. So my conjecture predicts
that $CH^*BPGL(3)$ is also not generated by Chern classes
(as happened in the case of the finite group $S_6$; see section
\ref{surjcases}).

Finally, we can begin to compute the Chow ring of the exceptional group
$G_2$. We use that $G_2$ is the group of automorphisms of
a general skew-symmetric cubic form on $\C^7$ \cite{FH},
p.~357. Equivalently, $GL(7)/G_2$ is a Zariski open subset
of $\Lambda^3(\C^7)=\C^{35}$, and so it has trivial Chow groups.
By Proposition \ref{gens}, it follows    
that the Chow ring of $BG_2$ is generated by the Chern classes
$c_1,\ldots,c_7$ of the standard representation $G_2\inj GL(7)$.


\end{document}